\documentclass[11pt]{article}
\usepackage{amsmath}
\usepackage{amsfonts}
\usepackage{amssymb}
\usepackage{float}
\usepackage{graphicx}
\usepackage{epstopdf}
\usepackage{hyperref}
\usepackage{multirow}
\usepackage[misc]{ifsym}

\providecommand{\U}[1]{\protect\rule{.1in}{.1in}}
\providecommand{\U}[1]{\protect\rule{.1in}{.1in}}
\providecommand{\U}[1]{\protect\rule{.1in}{.1in}}
\providecommand{\U}[1]{\protect\rule{.1in}{.1in}}
\newtheorem{theorem}{Theorem}

\newtheorem{algorithm}[theorem]{Algorithm}

\newtheorem{condition}[theorem]{Condition}

\numberwithin{equation}{section} \numberwithin{theorem}{section}
\newtheorem{definition}[theorem]{Definition}
\newtheorem{example}[theorem]{Example}

\newtheorem{lemma}[theorem]{Lemma}

\newtheorem{remark}[theorem]{Remark}

\begin{document}

\title{Bounded perturbation resilience of extragradient-type methods and their applications}

\author{ \textsc{Q.-L. Dong$^{a}$, A. Gibali$^{b}$, D. Jiang$^{a}$, Y. Tang$^{c}$}
\\
%EndAName
{\small ${^a}$Tianjin Key Laboratory for Advanced Signal Processing, College
of Science, }\\
{\small Civil Aviation University of China, Tianjin 300300, China, dongql@lsec.cc.ac.cn}\\
{\small $^b$(\Letter)\hspace{.01in}Department of Mathematics, ORT Braude College,}\\
{\small  2161002 Karmiel, Israel, avivg@braude.ac.il}\\
\small $^{c}$Department of Mathematics, NanChang University,\\
\small  Nanchang 330031, P.R. China, hhaaoo1331@163.com\\
}
\date{August 23, 2017. Revised: October 18, 2017. Accepted for publication in \textit{Journal of Inequalities and Applications}}
\maketitle

{} \textbf{Abstract.} In this paper we study the bounded perturbation
resilience of the extragradient and the subgradient extragradient methods
for solving variational inequality (VI) problem in real Hilbert spaces. This is an important property of algorithms
which guarantees the convergence of the scheme under summable errors, meaning that an inexact version of the methods can also be considered. Moreover, once an algorithm is proved to be bounded perturbation
resilience, superiorizion can be used, and this allows flexibility in choosing the bounded perturbations
in order to obtain a superior solution, as well explained in the paper. We also discuss some inertial extragradient methods.
Under mild and standard assumptions of monotonicity and Lipschitz continuity
of the VI's associated mapping, convergence of the perturbed extragradient
and subgradient extragradient methods is proved. In addition we show that
the perturbed algorithms converges at the rate of $O(1/t)$. Numerical
illustrations are given to demonstrate the performances of the algorithms. \\%
[6pt]

\section{Introduction}

In this paper we are concerned with the variational inequality (VI) problem
of finding a point $x^{\ast }$ such that%
\begin{equation}
\langle F(x^{\ast }),x-x^{\ast }\rangle \geq 0,\ \text{for all }x\in C,
\label{a}
\end{equation}%
where $C\subseteq \mathcal{H}$ is nonempty, closed and convex set in a real Hilbert
space $\mathcal{H}$, $\langle \cdot ,\cdot \rangle $ denotes the inner product in $\mathcal{H}$,
and $F:\mathcal{H}\rightarrow \mathcal{H}$ is a given mapping. This problem is a fundamental
problem in optimization theory and captures various applications, such as partial differential equations, optimal control, and mathematical programming; for theory and application of VIs and related problems the reader is referred for example to the works of Ceng et. al. \cite{clyy17}, Zegeye et al. \cite{zsy15}, the papers of Yao et. al. \cite{yly17, ynlk12, yply16} and the many references therein.

Many algorithms for solving the VI (\ref{a}) are projection algorithms that
employ projections onto the feasible set $C$ of the VI (\ref{a}), or onto
some related set, in order to reach iteratively a solution. Korpelevich \cite%
{Korpelevich} and Antipin \cite{Antipin} proposed an algorithm for solving
(\ref{a}), known as the \textit{extragradient method}, see also Facchinei
and Pang \cite[Chapter 12]{PF03}. In each iteration of the algorithm, in
order to get the next iterate $x^{k+1}$, two orthogonal projections onto $C$
are calculated, according to the following iterative step. Given the current
iterate $x^{k},$ calculate%
\begin{equation}
\left\{
\begin{array}{l}
y^{k}=P_{C}(x^{k}-\gamma _{k}F(x^{k}))\medskip \\
x^{k+1}=P_{C}(x^{k}-\gamma _{k}F(y^{k}))%
\end{array}%
\right.  \label{EG-1}
\end{equation}%
where $\gamma _{k}\in (0,1/L)$, and $L$ is the Lipschitz constant of $F$, or
$\gamma _{k}$ is updated by the following adaptive procedure
\begin{equation}
\gamma _{k}\Vert F(x^{k})-F(y^{k})\Vert \leq \mu \Vert x^{k}-y^{k}\Vert
,\quad \mu \in (0,1).  \label{a2}
\end{equation}%
In the extragradient method there is the need to calculate twice
the orthogonal projection onto $C$ in each iteration. In case that the set $%
C $ is simple enough so that projections onto it can be easily computed,
then this method is particularly useful; but if $C$ is a general closed and
convex set, a minimal distance problem has to be solved (twice) in order to
obtain the next iterate. This might seriously affect the efficiency of the
extragradient method. Hence, Censor et al in \cite{CGR, CGR1, CGR2}
presented a method called the \textit{subgradient extragradient method} in
which the second projection (\ref{EG-1}) onto $C$ is replaced by a specific
subgradient projection which can be easily calculated. The iterative step
has the following form.

\begin{equation}
\left\{
\begin{array}{l}
y^{k}=P_{C}(x^{k}-\gamma F(x^{k}))\medskip \\
x^{k+1}=P_{T_{k}}(x^{k}-\gamma F(y^{k}))%
\end{array}%
\right.  \label{alg:cgr}
\end{equation}%
where $T_{k}$ is the set defined as%
\begin{equation}
T_{k}:=\{w\in \mathcal{H}\mid \left\langle \left( x^{k}-\gamma F(x^{k})\right)
-y^{k},w-y^{k}\right\rangle \leq 0\},  \label{eq:T_K}
\end{equation}
and $\gamma\in(0,1/L).$

In this manuscript we prove that the above methods, the extragradient and the
subgradient extragradient methods are bounded perturbation resilient and the perturbed
methods have convergence rate of $O(1/t)$. This means that that will show that an inexact version of the algorithms, such that allows to incorporate summable errors also converge to a solution of the VI (\ref{a}) and moreover, their superiorized version can be introduced, by choosing the perturbations and in order to obtain a superior solution with respect to some new objective function, for example by choosing the norm, we can obtain a solution to the VI (\ref{a}) which is closer to the origin.

Our paper is organized as follows. In Section \ref{sec:pre} we present the
preliminaries. In Section \ref{sec:extragrad_per} we study the convergence
of the extragradient method with outer perturbations. Later in Section \ref%
{sec:bpr_eg} the bounded perturbation resilience of the extragradient method
is presented as well as the construction of the inertial extragradient
methods.

In the same spirit of the previous sections, in Section \ref{sec:sem} we
study the convergence of the subgradient extragradient method with outer
perturbations, show its bounded perturbation resilience and the construction
of the inertial subgradient extragradient methods. Finally, in Section \ref%
{sec:num} we present numerical examples in signal processing which demonstrate the performances
of the perturbed algorithms.

\section{Preliminaries\label{sec:pre}}

Let $\mathcal{H}$ be a real Hilbert space with inner product $\langle \cdot
,\cdot \rangle $ and the induced norm $\Vert \cdot \Vert $, and let $D$ be a
nonempty, closed and convex subset of $\mathcal{H}$. We write $%
x^{k}\rightharpoonup x$ to indicate that the sequence $\left\{ x^{k}\right\}
_{k=0}^{\infty }$ converges weakly to $x$ and $x^{k}\rightarrow x$ to
indicate that the sequence $\left\{ x^{k}\right\} _{k=0}^{\infty }$
converges strongly to $x.$ Given a sequence $\left\{ x^{k}\right\}
_{k=0}^{\infty }$, denote by $\omega _{w}(x^{k})$ its weak $\omega $-limit
set, that is, any $x\in \omega _{w}(x^{k})$ such that there exsists a subsequence $%
\left\{ x^{k_{j}}\right\} _{j=0}^{\infty }$of $\left\{ x^{k}\right\}
_{k=0}^{\infty }$ which converges weakly to $x$.

For each point $x\in \mathcal{H},$\ there exists a unique nearest point in $%
D $, denoted by $P_{D}(x)$. That is,%
\begin{equation}
\left\Vert x-P_{D}\left( x\right) \right\Vert \leq \left\Vert x-y\right\Vert
\text{ for all }y\in D.
\end{equation}%
The mapping $P_{D}:\mathcal{H}\rightarrow D$ is called the metric projection
of $\mathcal{H}$ onto $D$. It is well known that $P_{D}$ is a \textit{%
nonexpansive mapping} of $\mathcal{H}$ onto $D$, i.e., and even \textit{%
firmly nonexpansive mapping. This is captured in the next lemma.}

\begin{lemma}
\label{lem22} For any $x,y\in \mathcal{H}$ and $z\in D$, it holds

\begin{itemize}
\item $\Vert P_{D}(x)-P_{D}(y)\Vert ^{2}\leq \Vert x-y\Vert ;\medskip $

\item $\Vert P_{D}(x)-z\Vert ^{2}\leq \Vert x-z\Vert ^{2}-\Vert
P_{D}(x)-x\Vert ^{2}$;
\end{itemize}
\end{lemma}

The characterization of the metric projection $P_{D}$ \cite[Section 3]%
{Goebel+Reich}, is given by the following two properties in this lemma.

\begin{lemma}
\label{lem21} Given $x\in \mathcal{H}$ and $z\in D$. Then $z=P_{D}\left(
x\right) $ if and only if
\begin{equation}
P_{D}(x)\in D
\end{equation}%
and%
\begin{equation}
\left\langle x-P_{D}\left( x\right) ,P_{D}\left( x\right) -y\right\rangle
\geq 0\text{ for all }x\in \mathcal{H},\text{ }y\in D.
\end{equation}
\end{lemma}

\begin{definition}
The \texttt{normal cone} of $D$ at $v\in D$, denote by $N_{D}\left( v\right)
$ is defined as%
\begin{equation}
N_{D}\left( v\right) :=\{d\in \mathcal{H}\mid \left\langle
d,y-v\right\rangle \leq 0\text{ for all }y\in D\}.  \label{eq:normal_cone}
\end{equation}
\end{definition}

\begin{definition}
Let $B:\mathcal{H}\rightrightarrows 2^{\mathcal{H}}\mathcal{\ }$be a
point-to-set operator defined on a real Hilbert space $\mathcal{H}$. The
operator $B$ is called a \texttt{maximal monotone operator} if $B$ is
\texttt{monotone}, i.e.,%
\begin{equation}
\left\langle u-v,x-y\right\rangle \geq 0\text{ for all }u\in B(x)\text{ and }%
v\in B(y),
\end{equation}%
and the graph $G(B)$ of $B,$%
\begin{equation}
G(B):=\left\{ \left( x,u\right) \in \mathcal{H}\times \mathcal{H}\mid u\in
B(x)\right\} ,
\end{equation}%
is not properly contained in the graph of any other monotone operator.
\end{definition}

Based on Rockafellar (\cite[Theorem 3]{Rockafellar}), a monotone mapping $B$ is
maximal if and only if, for any $\left( x,u\right) \in\mathcal{H}\times%
\mathcal{H},$ if $\left\langle u-v,x-y\right\rangle \geq0$ for all $\left(
v,y\right) \in G(B)$, then it follows that $u\in B(x).$

\begin{definition}
The \texttt{subdifferential set} of a convex function $c$ at a point $x$ is defined as%
\begin{equation}
\partial c(x):=\{\xi \in \mathcal{H}\mid c(y)\geq c(x)+\langle \xi
,y-x\rangle \text{ for all }y\in \mathcal{H}\}.
\end{equation}
\end{definition}

For $z\in \mathcal{H},$ take any $\xi \in \partial c(z)$ and define%
\begin{equation}
T\left( z\right) :=\left\{ w\in \mathcal{H}\mid c(z)+\langle \xi ,w-z\rangle
\leq 0\right\} .  \label{eq:HyperS}
\end{equation}%
This is a half-space the bounding hyperplane of which separates the set $D$
from the point $z$ if $\xi \neq 0$; otherwise $T(z)=\mathcal{H};$ see, e.g.,
\cite[Lemma 7.3]{BB96}.

\begin{lemma}
\cite{CAM} \label{lem24} Let $D$ be a nonempty, closed and convex subset of
a Hilbert space $\mathcal{H}$. Let $\{x^{k}\}_{k=0}^{\infty }$ be a bounded
sequence which satisfies the following properties:

\begin{itemize}
\item every limit point of $\{x^{k}\}_{k=0}^{\infty }$ lies in $D$;

\item $\lim_{n\rightarrow \infty }\Vert x^{k}-x\Vert $ exists for every $%
x\in D$.
\end{itemize}
Then $\{x^{k}\}_{k=0}^{\infty }$ converges to a point in $D$.
\end{lemma}

\begin{lemma}
\label{lemineq} Assume that $\{a_{k}\}_{k=0}^{\infty }$ is a sequence of
nonnegative real numbers such that
\begin{equation}
a_{k+1}\leq (1+\gamma _{k})a_{k}+\delta _{k},\text{ }\quad \forall k\geq 0,
\end{equation}%
where the nonnegative sequences $\{\gamma _{k}\}_{k=0}^{\infty }$ and $%
\{\delta _{k}\}_{k=0}^{\infty }$ satisfy $\sum_{k=0}^{\infty }\gamma
_{k}<+\infty $ and $\sum_{k=0}^{\infty }\delta _{k}<+\infty $, respectively.
Then $\lim_{k\rightarrow \infty }a_{k}$ exists.
\end{lemma}

\section{The extragradient method with outer perturbations\label%
{sec:extragrad_per}}

In order to discuss the convergence of the extragradient method with outer
perturbations we make the following assumptions.\\

\begin{condition}
\label{con:Condition 1.1} The solution set of (\ref{a}), denoted by $SOL(C,F)$,
is nonempty.
\end{condition}

\begin{condition}
\label{con:Condition 1.2} The mapping $F$ is \textit{monotone} on $C$, i.e.,
\begin{equation}
\langle F(x)-F(y),x-y\rangle \geq 0,\quad \forall x,y\in C,  \label{a3}
\end{equation}
\end{condition}

\begin{condition}
\label{con:Condition 1.3} The mapping $F$ is\textit{\ Lipschitz continuous} on $%
C$ with the Lipschitz constant $L>0$, i.e.,
\begin{equation}
\Vert F(x)-F(y)\Vert \leq L\Vert x-y\Vert ,\quad \forall x,y\in C.
\label{a4}
\end{equation}
\end{condition}

Observe that while Censor et al in \cite[Theorem 3.1]{CGR1} showed the
weak convergence of the extragradient method (\ref{EG-1}) in Hilbert spaces
for a fixed step size $\gamma _{k}=\gamma \in (0,1/L)$, this can be easily
improved in case that the adaptive rule (\ref{a2}) is used. The next theorem shows this and its proof can easily be derived by following similar lines of the proof of \cite[Theorem 3.1]{CGR1}.

\begin{theorem}
Assume that Conditions \ref{con:Condition 1.1}--\ref{con:Condition 1.3} hold. Then any sequence $\{x^{k}\}_{k=0}^{%
\infty }$ generated by the extragradient method (\ref{EG-1}) with the adaptive
rule (\ref{a2}) weakly converges to a solution of the variational inequality
(\ref{a}).
\end{theorem}

Denote $e_{i}^{k}:=e_{i}(x^{k})$, $i=1,2$. The sequences of perturbations $%
\{e_{i}^{k}\}_{k=0}^{\infty }$, $i=1,2$, are assumed to be summable, i.e.,
\begin{equation}
\sum_{k=0}^{\infty }\Vert e_{i}^{k}\Vert <+\infty ,\quad i=1,2.
 \label{BP}
\end{equation}

Now we consider the extragradient method with outer perturbations.

\begin{algorithm}
\label{alg:em_per}$\left. {}\right. $\textbf{The extragradient method with
outer perturbations}

\textbf{Step 0}: Select a starting point $x^{0}\in C$
and set $k=0$.

\textbf{Step 1}: Given the current iterate $x^{k},$ compute%
\begin{equation}
\label{sp1}
y^{k}=P_{C}(x^{k}-\gamma _{k}F(x^{k})+e_{1}(x^{k})),
\end{equation}%
where $\gamma _{k}=\sigma \rho ^{m_{k}}$, $\sigma >0,$ $\rho \in (0,1)$ and $%
m_{k}$ is the smallest nonnegative integer such that (see \cite{Khobotov})
\begin{equation}
\gamma _{k}\Vert F(x^{k})-F(y^{k})\Vert \leq \mu \Vert x^{k}-y^{k}\Vert
,\quad \mu \in (0,1).
\label{ad4}
\end{equation}%
Calculate the next iterate%
\begin{equation}
x^{k+1}=P_{C}(x^{k}-\gamma _{k}F(y^{k})+e_{2}(x^{k})).
 \label{ai}
\end{equation}

\textbf{Step 2}:\textbf{\ }If\textbf{\ }$x^{k}=y^{k},$ then stop. Otherwise,
set $k\leftarrow (k+1)$ and return to \textbf{Step 1.}\bigskip
\end{algorithm}

\subsection{Convergence analysis}

\begin{lemma}
\cite{Yang} \label{lem31} The Armijo-like search rule (\ref{ad4}) is well
defined. Besides, $\underline{\gamma}\leq\gamma_k\leq \sigma$, where $%
\underline{\gamma}=\min\{\sigma,\frac{\mu\rho}{L}\}.$
\end{lemma}

\begin{theorem}
\label{Th31} Assume that Conditions \ref{con:Condition 1.1}--\ref{con:Condition 1.3} hold. Then the sequence $%
\{x^{k}\}_{k=0}^{\infty }$ generated by Algorithm \ref{alg:em_per} converges
weakly to a solution of the variational inequality (\ref{a}).
\end{theorem}

\textit{Proof.} Take $x^*\in SOL(C,F).$ From (\ref{ai}) and Lemma \ref%
{lem22}(ii), we have
\begin{equation}  \label{cc1}
\aligned \|x^{k+1}-x^*\|^2&\leq\|x^k-\gamma_k F(y^k)+e_2^k-x^*\|^2-
\|x^k-\gamma_k F(y^k)+e_2^k-x^{k+1}\|^2 \\
&=\|x^k-x^*\|^2-\|x^k-x^{k+1}\|^2+2\gamma_k\langle F(y^k),x^*-x^{k+1}\rangle\\
&\quad-2\langle e_2^k,x^*-x^{k+1}\rangle. \\
\endaligned
\end{equation}
From Cauchy-Schwartz inequality and the mean value inequality, it follows
\begin{equation}  \label{cc2}
\aligned -2\langle e_2^k,x^*-x^{k+1}\rangle &\leq2\|e_2^k\|\|x^{k+1}-x^*\| \\
&\leq\|e_2^k\|+\|e_2^k\|\|x^{k+1}-x^*\|^2. \endaligned
\end{equation}
Using $x^*\in SOL(C,F)$ and the monotone property of $F$, we have $\langle
y^k-x^*,F(y^k)\rangle\geq0$ and consequently get
\begin{equation}  \label{g1}
2\gamma_k\langle F(y^k),x^*-x^{k+1}\rangle\leq2\gamma_k\langle
F(y^k),y^k-x^{k+1}\rangle.
\end{equation}
Thus, we have
\begin{equation}  \label{c3}
\aligned -\|x^k-x^{k+1}\|^2&+2\gamma_k\langle F(y^k),x^*-x^{k+1}\rangle \\
&\leq-\|x^k-x^{k+1}\|^2+2\gamma_k\langle F(y^k),y^k-x^{k+1}\rangle \\
&=-\|x^k-y^k\|^2-\|y^k-x^{k+1}\|^2\\
&\quad+2\langle
x^k-\gamma_kF(y^k)-y^k,x^{k+1}-y^k\rangle, \endaligned
\end{equation}
where the equality comes from
\begin{equation}  \label{g1}
-\|x^k-x^{k+1}\|^2=-\|x^k-y^k\|^2-\|y^k-x^{k+1}\|^2-2\langle
x^k-y^k,y^k-x^{k+1}\rangle.
\end{equation}
Using $x^{k+1}\in C,$ the definition of $y^k$ and Lemma \ref{lem21}, we have
\begin{equation}  \label{g3}
\langle y^k-x^k+\gamma_kF(x^k)-e_1^k,x^{k+1}-y^k\rangle\geq0.
\end{equation}
So, we obtain
\begin{equation}  \label{c4}
\aligned 2\langle x^k-&\gamma_kF(y^k)-y^k,x^{k+1}-y^k\rangle \\
&\leq2\gamma_k\langle F(x^k)-F(y^k),x^{k+1}-y^k\rangle-2\langle
e_1^k,x^{k+1}-y^k\rangle \\
&\leq2\gamma_k\| F(x^k)-F(y^k)\|\|x^{k+1}-y^k\|+2\| e_1^k\|\|x^{k+1}-y^k\| \\
&\leq2\mu\| x^k-y^k\|\|x^{k+1}-y^k\|+\| e_1^k\|+\| e_1^k\|\|x^{k+1}-y^k\|^2
\\
&\leq\mu\| x^k-y^k\|^2+\mu\|x^{k+1}-y^k\|^2+\| e_1^k\|+\|
e_1^k\|\|x^{k+1}-y^k\|^2 \\
&=\mu\| x^k-y^k\|^2+(\mu+\| e_1^k\|)\|x^{k+1}-y^k\|^2+\| e_1^k\|. \endaligned
\end{equation}
From (\ref{BP}), it follows
\begin{equation}  \label{g4}
\lim_{k\rightarrow\infty}\|e_i^k\|=0,\quad i=1,2.
\end{equation}
Therefor, we assume $\|e_1^k\|\in[0,1- \mu-\nu)$ and $\|e_2^k\|\in[0,1/2)$, $%
k\geq0$, where $\nu\in (0,1-\mu)$. So, using (\ref{c4}), we get
\begin{equation}  \label{c4a}
2\langle x^k-\gamma_kF(y^k)-y^k,x^{k+1}-y^k\rangle\leq\mu\|
x^k-y^k\|^2+(1-\nu)\|x^{k+1}-y^k\|^2+\| e_1^k\|.
\end{equation}
Combining (\ref{cc1})-(\ref{c3}) and (\ref{c4a}), we obtain
\begin{equation}  \label{c5a}
\aligned
\|x^{k+1}-x^*\|^2 &\leq\|x^k-x^*\|^2-\left(1-\mu\right)\|
x^k-y^k\|^2-\nu\|x^{k+1}-y^k\|^2\\
&\quad+\|e^k\|+\|e_2^k\|\|x^{k+1}-x^*\|^2,
\endaligned
\end{equation}
where
\begin{equation}
\label{Th31-1}
\|e^k\|:=\|e_1^k\|+\|e_2^k\|.
\end{equation}
From (\ref{c5a}), it follows
\begin{equation}  \label{c5}
\aligned \|x^{k+1}-x^*\|^2&\leq \frac1{1-\|e_2^k\|}\|x^k-x^*\|^2 -\frac{1-\mu%
}{1-\|e_2^k\|}\| x^k-y^k\|^2 \\
&\quad-\frac{\nu}{1-\|e_2^k\|}\|x^{k+1}-y^k\|^2+\frac{\|e^k\|}{1-\|e_2^k\|}. %
\endaligned
\end{equation}
Since $\|e_2^k\|\in[0,1/2),$ $k\geq0$, we get
\begin{equation}  \label{g5}
1\leq\frac1{1-\|e_2^k\|}\leq1+2\|e_2^k\|<2.
\end{equation}
So, from (\ref{c5}), we have
\begin{equation}  \label{c6}
\aligned \|x^{k+1}-x^*\|^2&\leq (1+2\|e_2^k\|)\|x^k-x^*\|^2-(1-\mu)\|
x^k-y^k\|^2\\
&\quad-\nu\|x^{k+1}-y^k\|^2+2\|e^k\| \\
&\leq (1+2\|e_2^k\|)\|x^k-x^*\|^2+2\|e^k\|. \endaligned
\end{equation}
Using (\ref{BP}) and Lemma \ref{lemineq}, we get the existence of $%
\lim_{k\rightarrow\infty}\|x^k-x^*\|^2$ and then the boundedness of $%
\{x^k\}_{k=0}^\infty. $ From (\ref{c6}), it follows
\begin{equation}  \label{g6}
(1-\mu)\|
x^k-y^k\|^2+\nu\|x^{k+1}-y^k\|^2\leq(1+2\|e_2^k\|)\|x^k-x^*\|^2-%
\|x^{k+1}-x^*\|^2+2\|e^k\|,
\end{equation}
which means that
\begin{equation}  \label{g7}
\sum_{k=0}^\infty\| x^k-y^k\|^2 <+\infty,\quad \hbox{and}\quad
\sum_{k=0}^\infty\| x^{k+1}-y^k\|^2 <+\infty.
\end{equation}
Thus, we obtain
\begin{equation}  \label{c7}
\lim_{k\rightarrow\infty}\|x^k-y^k\|=0,\quad \hbox{and}\quad
\lim_{k\rightarrow\infty}\| x^{k+1}-y^k\|=0,
\end{equation}
and consequently,
\begin{equation}  \label{cd7}
\lim_{k\rightarrow\infty}\|x^{k+1}-x^k\|=0.
\end{equation}

Now, we are to show $\omega_w(x^k)\subseteq SOL(C,F).$ Due to the
boundedness of $\{x^k\}_{k=0}^\infty$, it has at least one weak accumulation point. Let $%
\hat x\in \omega_w(x^k)$. Then there exists a subsequence $\{x^{k_i}\}_{i=0}^\infty$ of $%
\{x^k\}_{k=0}^\infty$ which converges weakly to $\hat x$. From (\ref{c7}), it follows
that $\{y^{k_i}\}_{i=0}^\infty$ also converges weakly to $\hat x.$

We will show that $\hat{x}$ is a solution of the variational inequality (\ref%
{a}). Let%
\begin{equation}
A(v)=\left\{
\begin{array}{cc}
F(v)+N_{C}\left( v\right) , & v\in C, \\
\emptyset , & v\notin C\text{,}%
\end{array}%
\right.  \label{cx10}
\end{equation}%
where $N_{C}(v)$ is the normal cone of $C$ at $v\in C$. It is known that $A$
is a maximal monotone operator and $A^{-1}(0)=SOL(C,F)$. If $(v,w)\in G(A)$,
then we have $w-F(v)\in N_{C}(v)$ since $w\in A(v)=F(v)+N_{C}(v)$. Thus it
follows that
\begin{equation}
\langle w-F(v),v-y\rangle \geq 0,\quad y\in C.  \label{g9}
\end{equation}%
Since $y^{k_{i}}\in C$, we have
\begin{equation}
\langle w-F(v),v-y^{k_{i}}\rangle \geq 0.  \label{g10}
\end{equation}

On the other hand, by the definition of $y^k$ and Lemma \ref{lem21}, it
follows that
\begin{equation}  \label{g11}
\langle x^k-\gamma_k F(x^k)+e_1^k-y^k, y^{k}-v\rangle\geq0,
\end{equation}
and consequently,
\begin{equation}  \label{g12}
\left\langle \frac{y^k-x^k}{\gamma_k}+F(x^k),v-y^k\right\rangle-\frac1{\gamma_k}\langle
e_1^k, v-y^k \rangle\geq0.
\end{equation}
Hence we have
\begin{equation}  \label{g13}
\aligned
&\langle w,v-y^{k_i}\rangle\\
&\geq\langle F(v),v-y^{k_i}\rangle \\
&\geq\langle F(v),v-y^{k_i}\rangle-\Big\langle \frac{y^{k_i}-x^{k_i}}{%
\gamma_{{k_i}}}+F(x^{k_i}),v-y^{k_i}\Big\rangle +\frac1{\gamma_{k_i}}\langle
e_1^{k_i}, v-y^{k_i} \rangle \\
&=\langle F(v)-F(y^{k_i}),v-y^{k_i}\rangle+\langle
F(y^{k_i})-F(x^{k_i}),v-y^{k_i}\rangle\\
&\quad-\Big\langle \frac{y^{k_i}-x^{k_i}}{%
\gamma_{k_i}},v-y^{k_i}\Big\rangle  +\frac1{\gamma_{k_i}}\langle e_1^{k_i}, v-y^{k_i} \rangle \\
&\geq\langle F(y^{k_i})-F(x^{k_i}),v-y^{k_i}\rangle-\Big\langle \frac{%
y^{k_i}-x^{k_i}}{\gamma_{k_i}},v-y^{k_i}\Big\rangle+\frac1{\gamma_{k_i}}%
\langle e_1^{k_i}, v-y^{k_i} \rangle, \\
\endaligned
\end{equation}
which implies
\begin{equation}  \label{g14}
\langle w,v-y^{k_i}\rangle\geq\langle F(y^{k_i})-F(x^{k_i}),v-y^{k_i}\rangle-%
\Big\langle \frac{y^{k_i}-x^{k_i}}{\gamma_{k_i}},v-y^{k_i}\Big\rangle%
+\frac1{\gamma_{k_i}}\langle e_1^{k_i}, v-y^{k_i} \rangle. \\
\end{equation}
Taking the limit as $i\rightarrow\infty$ in the above inequality and using
Lemma \ref{lem31}, we obtain
\begin{equation}  \label{g15}
\langle w,v-\hat x\rangle\geq0.
\end{equation}
Since $A$ is a maximal monotone operator, it follows that $\hat x\in
A^{-1}(0) = SOL(C,F)$. So, $\omega_w(x^k)\subseteq SOL(C,F).$

Since $\lim_{k\rightarrow\infty}\|x^{k}-x^*\|$ exists and $%
\omega_w(x^k)\subseteq SOL(C,F)$, using Lemma \ref{lem24}, we conclude that $%
\{x^k\}_{k=0}^\infty$ weakly converges a solution of the variational inequality (\ref{a}%
). This completes the proof. $\Box$

\subsection{Convergence rate}

\qquad Nemirovski \cite{Nemirovski} and Tseng \cite{Tseng} proved the $%
O(1/t) $ convergence rate of the extragradient method. In this subsection,
we present the convergence rate of Algorithm \ref{alg:em_per}.

\begin{theorem}
\label{Th32} Assume that Conditions \ref{con:Condition 1.1}--\ref{con:Condition 1.3} hold. Let the sequences $%
\{x^{k}\}_{k=0}^{\infty }$ and $\{y^{k}\}_{k=0}^{\infty }$ be generated by
Algorithm \ref{alg:em_per}. For any integer $t>0$, we have a $y_{t}\in C$
which satisfies
\begin{equation}
\aligned\langle F(x),y_{t}-x\rangle \leq \frac{1}{2 \Upsilon _{t}}%
(\Vert x-x^{0}\Vert ^{2}+M(x)),\quad \forall x\in C,\endaligned  \label{m-1}
\end{equation}%
where
\begin{equation}
\label{lem36-0}
y_{t}=\frac{1}{\Upsilon _{t}}\sum_{k=0}^{t}\gamma _{k}y^{k},\,\, \Upsilon
_{t}=\sum_{k=0}^{t}\gamma _{k}
\end{equation}
and
\begin{equation}
 M(x)=\sup_k\{\max\{\|x^{k+1}-y^k\|,\|x^{k+1}-x\|\}\}\sum_{k=0}^{\infty}\Vert e^{k}\Vert.
\end{equation}
\end{theorem}

\textit{Proof.} Take arbitrarily $x\in C.$ From Conditions \ref{con:Condition 1.2}
and \ref{con:Condition 1.3}, we have
\begin{equation}
\label{Th32-1}
\aligned
&-\|x^k-x^{k+1}\|^2+2\gamma_k\langle F(y^k),x-x^{k+1}\rangle\\
&=-\|x^k-x^{k+1}\|^2+2\gamma_k\big[\langle F(y^k)-F(x),x-y^k\rangle+\langle F(x),x-y^k\rangle\\
&\quad+\langle F(y^k),y^k-x^{k+1}\rangle\big]\\
&\leq-\|x^k-x^{k+1}\|^2+2\gamma_k\big[\langle F(x),x-y^k\rangle+\langle F(y^k),y^k-x^{k+1}\rangle\big]\\
&=-\|x^k-y^k\|^2-\|y^k-x^{k+1}\|^2+2\gamma_k\langle F(x),x-y^k\rangle\\
&\quad+2\langle x^k-\gamma_kF(y^k)-y^k,x^{k+1}-y^k\rangle.
\endaligned
\end{equation}
By (\ref{ai}) and Lemma \ref{lem21}, we get
\begin{equation}
\label{Th32-2}
\aligned
&2\langle x^k-\gamma_kF(y^k)-y^k,x^{k+1}-y^k\rangle\\
&=2\langle x^k-\gamma_kF(x^k)+e_1^k-y^k,x^{k+1}-y^k\rangle-2\langle e_1^k,x^{k+1}-y^k\rangle\\
&\quad+2\gamma_k\langle F(x^k)-F(y^k),x^{k+1}-y^k\rangle\\
&\leq-2\langle e_1^k,x^{k+1}-y^k\rangle+2\gamma_k\langle F(x^k)-F(y^k),x^{k+1}-y^k\rangle\\
&\leq2\|e_1^k\|\|x^{k+1}-y^k\|+2\mu\|x^k-y^k\|\|x^{k+1}-y^k\|\\
&\leq2\|e_1^k\|\|x^{k+1}-y^k\|+\mu^2\|x^k-y^k\|^2+\|x^{k+1}-y^k\|^2.\\
\endaligned
\end{equation}
Identifying $x^*$ with $x$ in (\ref{cc1}) and (\ref{cc2}), and combining (%
\ref{Th32-1}) and (\ref{Th32-2}), we get
\begin{equation}
\label{Th32-5}
\aligned
&\|x^{k+1}-x\|^2\\
&\leq\|x^k-x\|^2+2\|e_1^k\|\|x^{k+1}-y^k\|-(1-\mu^2)\|x^k-y^k\|^2\\
&\quad+2\|e_2^k\|\|x^{k+1}-x\|+2\gamma_k\langle F(x),x-y^k\rangle\\
&\leq\|x^k-x\|^2+2\|e_1^k\|\|x^{k+1}-y^k\|+2\|e_2^k\|\|x^{k+1}-x\|\\
&\quad+2\gamma_k\langle F(x),x-y^k\rangle.\\
\endaligned
\end{equation}
Thus, we have
\begin{equation}
\label{m2}
\aligned
&\gamma_k \langle
F(x),y^k-x\rangle\\
&\leq\frac12(\|x^k-x\|^2-\|x^{k+1}-x\|^2)+\|e_1^k\|\|x^{k+1}-y^k\|+\|e_2^k\|\|x^{k+1}-x\| \\
&\leq\frac12(\|x^k-x\|^2-\|x^{k+1}-x\|^2)+M^{\prime }(x)\|e^k\| \\
\endaligned
\end{equation}
where $M^{\prime }(x)=\sup_k\{\max\{\|x^{k+1}-y^k\|,\|x^{k+1}-x\|\}\}<+\infty$. Summing the
inequality (\ref{m2}) over $k = 0,\ldots, t$, we obtain
\begin{equation}  \label{m3}
\aligned \left\langle
F(x),\sum_{k=0}^t\gamma_ky^k-\left(\sum_{k=0}^t\gamma_k\right)x\right\rangle
&\leq\frac12\|x^0-x\|^2+\frac{M^{\prime
}(x)}2\sum_{k=0}^t\|e^k\| \\
&=\frac12\|x^0-x\|^2+\frac12M(x). \\
\endaligned
\end{equation}
Using the notations of $\Upsilon_t$ and $y^t$ in the above inequality, we
derive
\begin{equation}  \label{m4}
\langle F(x),y_t-x\rangle\leq\frac{1}{2\Upsilon_t}(\|x-x^0\|^2+M(x)),%
\quad \forall x\in C.
\end{equation}
 The proof is complete. $\Box$

\begin{remark}
\textrm{\ From Lemma \ref{lem31}, it follows
\begin{equation}  \label{m5}
\Upsilon_t\geq(t+1)\underline{\gamma},
\end{equation}
thus Algorithm \ref{alg:em_per} has $O(1/t)$
convergence rate. In fact, for any bounded subset $D \subset C$ and given accuracy $\epsilon > 0$, our
algorithm achieves
\begin{equation}
\label{remark39-1}
\langle F(x),y_t-x\rangle\leq\epsilon,\quad \forall x\in D
\end{equation}
in at most
\begin{equation}
\label{remark39-2}
t=\left[\frac{m}{2\underline{\gamma}\epsilon}\right]
\end{equation}
iterations, where $y_t$ is defined by (\ref{lem36-0}) and $m=\sup\{\|x-x^0\|^2+M(x)|x\in D\}.$
}
\end{remark}

\section{The bounded perturbation resilience of the extragradient method
\label{sec:bpr_eg}}

In this section, we prove the bounded perturbation resilience (BPR) of
the extragradient method. This property is fundamental for the application
of the superiorization methodology (SM) to them.

The \textit{superiorization methodology} first appeared in Butnariu et al. in \cite{bdhk07}, without mentioning specifically the words superiorization and perturbation resilience. Some of the results in \cite{bdhk07} are based on earlier results of Butnariu, Reich and Zaslavski \cite{brz06, brz07, brz08}. For the state of current research on superiorization, visit the webpage:
\textquotedblleft Superiorization and Perturbation Resilience of Algorithms: A Bibliography compiled and continuously updated by Yair Censor\textquotedblright\  at: \href{http://math.haifa.ac.il/yair/bib-superiorization-censor.html}{http://math.haifa.ac.il/yair/bib-superiorization-censor.html}
and in particular see \cite[Section 3]{Censor14} and \cite[Appendix]{Censor17}.

Originally, the superiorization methodology is intended for constrained minimization (CM) problems of the form:
\begin{equation}
\min\,\{\phi (x)\,|\,x\in \Psi \}  \label{e1}
\end{equation}%
where $\phi :H\rightarrow \mathbb{R}$ is an objective function and $\Psi
\subseteq H$ is the solution set another problem. Here, we assume $\Psi \neq
\emptyset $ throughout this paper. Assume that the set $\Psi $ is a closed
convex subsets of a Hilbert space $H$, the minimization problem (\ref{e1}) becomes a standard CM
problem. Here we are interested in the case wherein $\Psi $ is the solution
set of another CM of the form:
\begin{equation}
\min\,\,\{f(x)\,|\,x\in \Omega \}  \label{e2}
\end{equation}%
i.e., we wish to look at
\begin{equation}
\Psi :=\{x^{\ast }\in \Omega \,|\,f(x^{\ast })\leq f(x)\,|\,\hbox{for all}%
\,\,x\in \Omega \}  \label{e3}
\end{equation}%
provided that $\Psi $ is nonempty. If $f$ is differentiable and let $%
F=\nabla f$, then the CM (\ref{e2}) equals to the following variational
inequality: to find a point $x^{\ast }\in C$ such that
\begin{equation}
\langle F(x^{\ast }),x-x^{\ast }\rangle \geq 0,\quad \forall x\in C.
\label{e4}
\end{equation}

The superiorization methodology (SM) strives not to solve (\ref{e1}) but
rather the task is to find a point in $\Psi$ which is superior, i.e., has a
lower, but not necessarily minimal, value of the objective function $\phi$.
This is done in the SM by first investigating the bounded perturbation
resilience of an algorithm designed to solve (\ref{e2}) and then proactively
using such permitted perturbations in order to steer the iterates of such an
algorithm toward lower values of the $\phi$ objective function while not
loosing the overall convergence to a point in $\Psi$.

In this paper, we do not investigate superiorization of the extragradient
method. We prepare for such an application by proving the bounded
perturbation resilience that is needed in order to do superiorization.

\begin{algorithm}
\textbf{The Basic Algorithm}

\textbf{Initialization:} $x^0 \in \Theta$ is arbitrary;

\textbf{Iterative Step:} Given the current iterate vector $x^{k}$, calculate
the next iterate $x^{k+1}$ via
\begin{equation}  \label{BA1}
x^{k+1}=\mathbf{A}_\Psi(x^k).
\end{equation}
\end{algorithm}

The bounded perturbation resilience (henceforth abbreviated by BPR) of such
a basic algorithm is defined next.

\begin{definition}
\label{Def41} \cite{HGDC} An algorithmic operator $\mathbf{A}_\Psi : H\rightarrow \Theta$ is said to be \texttt{
bounded perturbations resilient} if the following is true. If Algorithm \ref{BA1}
generates sequences $\{x^k\}_{k=0}^\infty$ with $x^0\in \Theta, $ that
converge to points in $\Psi$, then any sequence $\{y^k\}_{k=0}^\infty$,
starting from any $y^0\in \Theta,$ generated by
\begin{equation}  \label{mBA}
y^{k+1} = \mathbf{A_\Psi}(y^k + \lambda_kv^k),\quad \hbox{for all}\,\, k
\geq0,
\end{equation}
also converges to a point in $\Psi$, provided that, (i) the sequence $%
\{v^k\}_{k=0}^\infty$ is bounded, and (ii) the scalars $\{\lambda_k\}_{k=0}^%
\infty$ are such that $\lambda_k\geq 0$ for all $k \geq 0$, and $%
\sum_{k=0}^\infty \lambda_k<+\infty$, and (iii) $y^k+\lambda_kv^k\in \Theta$
for all $k \geq 0$.
\end{definition}

Definition \ref{Def41} is non-trivial only if $\Theta\neq \mathcal{H}$, in which the
condition (iii) is enforced in the superiorized version of the basic
algorithm, see step (xiv) in the ``Superiorized Version of Algorithm P" in (%
\cite{HGDC}, p. 5537) and step (14) in ``Superiorized Version of the ML-EM
Algorithm" in (\cite{GH}, Subsection II.B). This will be the case in the
present work.

Treating the extragradient method as the Basic Algorithm $\mathbf{A}_\Psi$,
our strategy is to first prove convergence of the iterative step (\ref{EG-1}) with bounded perturbations.
We show next how the convergence of this yields BPR according to Definition \ref{Def41}.

A superiorized version of any Basic Algorithm employs the perturbed version
of the Basic Algorithm as in (\ref{mBA}). A certificate to do so in the
superiorization method, see \cite{Dong2}, is gained by showing that the Basic
Algorithm is BPR. Therefore, proving the BPR of an algorithm is the first
step toward superiorizing it. This is done for the extragradient method in
the next subsection.

\subsection{The BPR of the extragradient method}

\qquad In this subsection, we investigate the bounded perturbation
resilience of the extragradient method whose iterative step is given by (\ref%
{EG-1}).

To this end, we treat the right-hand side of (\ref{EG-1}) as the algorithmic
operator $\mathbf{A}_\Psi$ of Definition \ref{Def41}, namely, we define for
all $k\geq0,$
\begin{equation}  \label{g16}
\mathbf{A}_\Psi(x^k)=P_C(x^k-\gamma_kF(P_C(x^k-\gamma_k F(x^k)))),
\end{equation}
and identify the solution set $\Psi$ with the solution set of the
variational inequality (\ref{a}) and identify the additional set $\Theta$
with $C$.

According to Definition \ref{Def41}, we need to show the convergence of the
sequence $\{x^{k}\}_{k=0}^{\infty }$ that, starting from any $x^{0}\in C$,
is generated by
\begin{equation}
x^{k+1}=P_{C}((x^{k}+\lambda _{k}v^{k})-\gamma _{k}F(P_{C}((x^{k}+\lambda
_{k}v^{k})-\gamma _{k}F(x^{k}+\lambda _{k}v^{k})))),  \label{g17}
\end{equation}%
which can be rewritten as

\begin{equation}
\left\{
\begin{array}{l}
y^{k}=P_{C}((x^{k}+\lambda _{k}v^{k})-\gamma _{k}F(x^{k}+\lambda
_{k}v^{k}))\medskip \\
x^{k+1}=P_{C}((x^{k}+\lambda _{k}v^{k})-\gamma _{k}F(y^{k}))%
\end{array}%
\right.  \label{mEG}
\end{equation}%

where $\gamma _{k}=\sigma \rho ^{m_{k}}$, $\sigma >0,$ $\rho \in (0,1)$ and $%
m_{k}$ is the smallest nonnegative integer such that
\begin{equation}
\gamma _{k}\Vert F(x^{k}+\lambda _{k}v^{k})-F(y^{k})\Vert \leq \mu (\Vert
x^{k}-y^{k}\Vert +\lambda _{k}\Vert v^{k}\Vert ),\quad \mu \in (0,1).
\label{sa1}
\end{equation}%
The sequences $\{v^{k}\}_{k=0}^{\infty }$ and $\{\lambda _{k}\}_{k=0}^{\infty }$ obey
the conditions (i) and (ii) in Definition \ref{Def41}, respectively, and
also (iii) in Definition \ref{Def41} is satisfied.\vskip0.1cm

The next theorem establishes the bounded perturbation resilience of the
extragradient method. The proof idea is to build a relationship between BPR
and the convergence of the iterative step (\ref{EG-1}).

\begin{theorem}
\label{Th41} Assume that Conditions \ref{con:Condition 1.1}-\ref{con:Condition 1.3} hold. Assume the
sequence $\{v^k\}_{k=0}^\infty$ is bounded, and the scalars $%
\{\lambda_k\}_{k=0}^\infty$ are such that $\lambda_k\geq 0$ for all $k \geq
0 $, and $\sum_{k=0}^\infty \lambda_k<+\infty$. Then the sequence $\{x^k\}_{k=0}^\infty$
generated by (\ref{mEG}) and (\ref{sa1}) converges weakly to a solution of
the variational inequality (\ref{a}).
\end{theorem}

\textit{Proof.} Take $x^*\in SOL(C,F).$ From $\sum_{k=0}^\infty
\lambda_k<+\infty$ and that $\{v^k\}_{k=0}^\infty$ is bounded, we have
\begin{equation}  \label{d1}
\sum_{k=0}^\infty \lambda_k\|v^k\|<+\infty,
\end{equation}
which means
\begin{equation}  \label{g18}
\lim_{k\rightarrow\infty}\lambda_k\|v^k\|=0.
\end{equation}
So, we assume $\lambda_k\|v^k\|\in [0,(1-\mu-\nu)/2)$, where $\nu\in
[0,1-\mu)$. Identifying $e_2^k$ with $\lambda_kv^k$ in (\ref{cc1}) and (\ref%
{cc2}) and using (\ref{c3}), we get

\begin{equation}  \label{d6}
\aligned \|x^{k+1}-x^*\|^2&=\|x^k-x^*\|^2+\lambda_k\|v^k\|+\lambda_k\|v^k\|
\|x^{k+1}-x^*\|^2-\|x^k-y^k\|^2 \\
&\quad-\|y^k-x^{k+1}\|^2+2\langle x^k-\gamma_k F(y^k)-y^k,x^{k+1}-y^k\rangle. \\
\endaligned
\end{equation}

From $x^{k+1}\in C,$ the definition of $y^{k}$ and Lemma \ref{lem21}, we
have
\begin{equation*}
\aligned \langle y^{k}-x^{k}-\lambda _{k}v^{k}+\gamma _{k}F(x^{k}+\lambda
_{k}v^{k}),x^{k+1}-y^{k}\rangle \geq 0. \\
\endaligned
\end{equation*}%
So, we obtain
\begin{equation}
\aligned
&2\langle x^{k}-\gamma _{k}F(y^{k})-y^{k},x^{k+1}-y^{k}\rangle
\label{d4} \\
&\leq 2\gamma _{k}\langle F(x^{k}+\lambda
_{k}v^{k})-F(y^{k}),x^{k+1}-y^{k}\rangle -2\lambda _{k}\langle
v^{k},x^{k+1}-y^{k}\rangle . \\
\endaligned
\end{equation}%
We have
\begin{equation}
\aligned
2\gamma _{k}\langle F(x^{k}&+\lambda
_{k}v^{k})-F(y^{k}),x^{k+1}-y^{k}\rangle  \label{d3} \\
&\leq 2\gamma _{k}\Vert F(x^{k}+\lambda _{k}v^{k})-F(y^{k})\Vert \Vert
x^{k+1}-y^{k}\Vert \\
&\leq 2\mu \Vert x^{k}+\lambda _{k}v^{k}-y^{k}\Vert \Vert x^{k+1}-y^{k}\Vert
\\
&\leq 2\mu \left( \Vert x^{k}-y^{k}\Vert +\lambda _{k}\Vert v^{k}\Vert
\right) \Vert x^{k+1}-y^{k}\Vert \\
&\leq 2\mu \Vert x^{k}-y^{k}\Vert \Vert x^{k+1}-y^{k}\Vert +2\mu \lambda
_{k}\Vert v^{k}\Vert \Vert x^{k+1}-y^{k}\Vert \\
&\leq \mu \Vert x^{k}-y^{k}\Vert ^{2}+(\mu +\lambda _{k}\Vert v^{k}\Vert
)\Vert x^{k+1}-y^{k}\Vert ^{2}+\mu ^{2}\lambda _{k}\Vert v^{k}\Vert . \\
\endaligned
\end{equation}%
Similarly with (\ref{cc2}), we can show
\begin{equation}  \label{d2}
\aligned-2\lambda _{k}\langle v^{k},x^{k+1}-y^{k}\rangle \leq \lambda
_{k}\Vert v^{k}\Vert +\lambda _{k}\Vert v^{k}\Vert \Vert x^{k+1}-y^{k}\Vert
^{2}. \endaligned
\end{equation}%
Combining (\ref{d4})-(\ref{d2}), we get
\begin{equation}
\aligned
2\langle x^{k}&-\gamma _{k}F(y^{k})-y^{k},x^{k+1}-y^{k}\rangle
\label{d7} \\
&\leq \mu \Vert x^{k}-y^{k}\Vert ^{2}+(\mu +2\lambda _{k}\Vert v^{k}\Vert
)\Vert x^{k+1}-y^{k}\Vert ^{2}+(1+\mu ^{2})\lambda _{k}\Vert v^{k}\Vert \\
&\leq \mu \Vert x^{k}-y^{k}\Vert ^{2}+(1-\nu )\Vert x^{k+1}-y^{k}\Vert
^{2}+2\lambda _{k}\Vert v^{k}\Vert ,\endaligned
\end{equation}%
where the last inequality comes from $\lambda _{k}\Vert v^{k}\Vert <(1-\mu
)/2$ and $\mu <1$. Substituting (\ref{d7}) into (\ref{d6}), we get
\begin{equation}
\aligned
\Vert x^{k+1}-x^{\ast }\Vert ^{2}&\leq \Vert x^{k}-x^{\ast }\Vert
^{2}-(1-\mu )\Vert x^{k}-y^{k}\Vert ^{2}-\nu \Vert x^{k+1}-y^{k}\Vert
^{2}+3\lambda _{k}\Vert v^{k}\Vert  \label{d8} \\
&\quad +\Vert x^{k+1}-x^{\ast }\Vert ^{2}.
\endaligned
\end{equation}%
Following the proof line of Theorem \ref{Th31}, we get $\{x^{k}\}_{k=0}^\infty$ weakly
converges to a solution of the variational equality (\ref{a}). $\Box $%
\newline

By using Theorems \ref{Th32} and \ref{Th41}, we obtain the convergence rate
of the extragradient method with BP.

\begin{theorem}
\label{Th42} Assume that Conditions \ref{con:Condition 1.1}-\ref{con:Condition 1.3} hold. Assume the
sequence $\{v^k\}_{k=0}^\infty$ is bounded, and the scalars $%
\{\lambda_k\}_{k=0}^\infty$ are such that $\lambda_k\geq 0$ for all $k \geq
0 $, and $\sum_{k=0}^\infty \lambda_k<+\infty$. Let the sequences $\{x^k\}_{k=0}^\infty$
and $\{y^k\}_{k=0}^\infty$ be generated by (\ref{mEG}) and (\ref{sa1}). For any integer $%
t > 0$, we have a $y_t\in C$ which satisfies
\begin{equation}  \label{m-1}
\aligned \langle F(x),y_t-x\rangle\leq\frac{1}{2\Upsilon_t}%
(\|x-x^0\|^2+M(x)),\quad \forall x\in C,\endaligned
\end{equation}
where
\begin{equation}
\label{Th42-1}
y_t=\frac{1}{\Upsilon_t}\sum_{k=0}^t\gamma_ky^k,\quad
\Upsilon_t=\sum_{k=0}^t\gamma_k,
\end{equation}
and
\begin{equation}
\label{Th42-2}
M(x)=\sup_k\{\max\{\|x^{k+1}-y^k\|,3\|x^{k+1}-x\|^2\}\}\sum_{k=0}^\infty\lambda_k\|v^k\|.
\end{equation}
\end{theorem}

\subsection{Construction of the inertial extragradient methods by BPR}

\qquad In this subsection, we construct two class of inertial extragradient
methods by using BPR, i.e., identifying the $e_i^k$, $k=1,2$ and $\lambda_k$%
, $v^k$ with special values.

Polyak \cite{Polyak1,Polyak2} first introduced the inertial-type algorithms
by using the heavy ball method of the second-order dynamical systems in
time. Since the inertial-type algorithms speed up the original algorithms
without the inertial effects, recently there are increasing interests in
studying inertial-type algorithms, (see, e.g. \cite{Alvarez,APR,BC,OBP,OBP}%
). The authors \cite{Dong} introduced an inertial extragradient method as
follows:
\begin{equation}\label{Alg:iegm}
\left\{
\aligned
&w^{k}=x^{k}+\alpha _{k}(x^{k}-x^{k-1}),  \\
&y^{k}=P_{C}(w^{k}-\gamma F(w^{k})), \\
&x^{k+1}=(1-\lambda _{k})w^{k}+\lambda _{k}P_{C}(w^{k}-\gamma F(y^{k}))
\endaligned
\right.
\end{equation}%
for each $k\geq 1$, where $\gamma \in (0,1/L),$ $\{\alpha _{k}\}$ is
nondecreasing with $\alpha _{1}=0$ and $0\leq \alpha _{k}\leq \alpha <1$ for
each $k\geq 1$ and $\lambda ,\sigma ,\delta >0$ are such that
\begin{equation}
\delta >\frac{\alpha \lbrack (1+\gamma L)^{2}\alpha (1+\alpha )+(1-\gamma
^{2}L^{2})\alpha \sigma +\sigma (1+\gamma L)^{2}]}{1-\gamma ^{2}L^{2}}
\label{delta}
\end{equation}%
and
\begin{equation*}
0<\lambda \leq \lambda _{k}\leq \frac{\delta (1-\gamma ^{2}L^{2})-\alpha
\lbrack (1+\gamma L)^{2}\alpha (1+\alpha )+(1-\gamma ^{2}L^{2})\alpha \sigma
+\sigma (1+\gamma L)^{2}]}{\delta \lbrack (1+\gamma L)^{2}\alpha (1+\alpha
)+(1-\gamma ^{2}L^{2})\alpha \sigma +\sigma (1+\gamma L)^{2}]},
\end{equation*}%
where $L$ is the Lipschitz constant of $F$.

Based on the iterative step (\ref{EG-1}), we construct the following inertial extragradient method:
\begin{equation}
\left\{
\begin{array}{ll}
y^{k}=P_{C}(x^{k}-\gamma _{k}F(x^{k})+\alpha _{k}^{(1)}(x^{k}-x^{k-1})), \\
x^{k+1}=P_{C}(x^{k}-\gamma _{k}F(y^{k})+\alpha _{k}^{(2)}(x^{k}-x^{k-1})),%
\end{array}%
\right.\label{Alg:I-iegm}
\end{equation}%
where
\begin{equation}
\alpha _{k}^{(i)}=\left\{
\begin{array}{ll}
\frac{\beta _{k}^{(i)}}{\left\Vert x^{k}-x^{k-1}\right\Vert }, & \text{if }%
\left\Vert x^{k}-x^{k-1}\right\Vert >1,\text{ }i=1,2 \\
\beta _{k}^{(i)}, & \text{if }\left\Vert x^{k}-x^{k-1}\right\Vert \leq 1.%
\end{array}%
\right.
\end{equation}

\begin{theorem}
\label{Th51} Assume that Conditions \ref{con:Condition 1.1}--\ref{con:Condition 1.3} hold. Assume that the
sequences $\{\beta_k^{(i)}\}_{k=0}^\infty$, $i=1,2$ satisfy $%
\sum_{k=1}^{\infty}\beta_k^{(i)}<\infty$, $i=1,2$. Then the sequence $%
\{x^k\}_{k=0}^\infty$ generated by the inertial extragradient method \ref{Alg:I-iegm}
converges weakly to a solution of the variational inequality $(\ref{a})$.
\end{theorem}

\textit{Proof.} Let $e_{i}^{k}=\beta _{k}^{(i)}v^{k}$, $i=1,2,$ where
\begin{equation}
v^{k}=\left\{
\begin{array}{ll}
\frac{x^{k}-x^{k-1}}{\left\Vert x^{k}-x^{k-1}\right\Vert }, & \text{if }%
\left\Vert x^{k}-x^{k-1}\right\Vert >1,\text{ }i=1,2 \\
x^{k}-x^{k-1}, & \text{if }\left\Vert x^{k}-x^{k-1}\right\Vert \leq 1.%
\end{array}%
\right.
\end{equation}%
It is obvious that $\Vert v^{k}\Vert \leq 1.$ So, it follows that $%
\{e_{i}^{k}\}$, $i=1,2$ satisfy (\ref{BP}) from the condition on $%
\{\beta _{k}^{(i)}\}.$ Using Theorem \ref{Th31}, we complete the proof. $\Box$

\begin{remark}
\textrm{From (\ref{cd7}), we have $\|x^k- x^{k-1}\|\leq1$ for big enough $k$%
, that is $\alpha_k^{(i)}=\beta_k^{(i)}.$ }
\end{remark}

Using the extragradient method with bounded perturbations (\ref{mEG}), we
construct the following inertial extragradient method:%
\begin{equation}
\left\{
\begin{array}{ll}
y^{k}=P_{C}(x^{k}+\alpha _{k}(x^{k}-x^{k-1})-\gamma _{k}F(x^{k}+\alpha
_{k}(x^{k}-x^{k-1}))), \\
x^{k+1}=P_{C}(x^{k}+\alpha _{k}(x^{k}-x^{k-1})-\gamma _{k}F(y^{k})),%
\end{array}%
\right.  \label{Alg:II-iegm}
\end{equation}%
where%
\begin{equation}
\alpha _{k}=\left\{
\begin{array}{ll}
\frac{\beta _{k}}{\left\Vert x^{k}-x^{k-1}\right\Vert }, & \text{if }%
\left\Vert x^{k}-x^{k-1}\right\Vert >1,\text{ }i=1,2 \\
\beta _{k}, & \text{if }\left\Vert x^{k}-x^{k-1}\right\Vert \leq 1.%
\end{array}%
\right.  \label{vk2a}
\end{equation}%
We extend Theorem \ref{Th41} to the convergence of the inertial
extragradient method \ref{Alg:II-iegm}.

\begin{theorem}
\label{Th52} Assume that Conditions \ref{con:Condition 1.1}--\ref{con:Condition 1.3} hold. Assume that the
sequence $\{\beta _{k}\}_{k=0}^{\infty }$ satisfies $\sum_{k=1}^{\infty
}\beta _{k}<\infty $. Then the sequence $\{x^{k}\}_{k=0}^{\infty }$
generated by the inertial extragradient method (\ref{Alg:II-iegm}) converges weakly to a
solution of the variational inequality (\ref{a}).
\end{theorem}

\begin{remark}
\textrm{The inertial parameter $\alpha_k$ in the inertial extragradient
method (\ref{Alg:I-iegm}) is bigger than that of inertial extragradient method
(\ref{Alg:II-iegm}). The inertial extragradient method (\ref{Alg:I-iegm}) becomes the inertial
extragradient method (\ref{Alg:II-iegm}) when $\lambda_k=1.$ }
\end{remark}

\section{The extension to the subgradient extragradient method\label{sec:sem}%
}

\qquad In this section, we generalize the results of extragradient method
proposed in the previous sections to the subgradient extragradient method.

Censor \textit{et al.} \cite{CGR} presented the subgradient extragradient
method (\ref{alg:cgr}). In their method the step size is fixed
 $\gamma \in (0,1/L)$, where $L$ is Lipschitz constant of $F$. So, in
order to determine the stepsize $\gamma _{k}$, one needs first calculate (or
estimate) $L$, which might be difficult or even impossible in general. So,
in order to overcome this, armijo-like search rule can be used:%
\begin{equation}
\gamma _{k}\Vert F(x^{k})-F(y^{k})\Vert \leq \mu \Vert x^{k}-y^{k}\Vert
,\quad \mu \in (0,1).  \label{as}
\end{equation}

To discuss the convergence of the subgradient extragradient method, we make
the following assumptions:

\begin{condition}\label{con:Condition 1.4}
The mapping $F$ is monotone on $\mathcal{H}$, i.e.,
\begin{equation}
\langle f(x)-f(y),x-y\rangle \geq 0,\quad \forall x,y\in \mathcal{H},
\end{equation}
\end{condition}

\begin{condition}\label{con:Condition 1.5}
%\textbf{Condition 1.5}
The mapping $F$ is Lipschitz continuous on $\mathcal{H%
}$ with the Lipschitz constant $L>0$, i.e.,
\begin{equation}
\Vert F(x)-F(y)\Vert \leq L\Vert x-y\Vert ,\quad \forall x,y\in \mathcal{H}.
\end{equation}
\end{condition}

As before,  Censor et al's subgradient
extragradient method (\cite[Theorem 3.1]{CGR1})  can be easily generalized by using some adaptive step rule, for example (%
\ref{as}). This result is captured in the next theorem.

\begin{theorem}
Assume that Conditions \ref{con:Condition 1.1}, \ref{con:Condition 1.4} and \ref{con:Condition 1.5} hold. Then the sequence $%
\{x^{k}\}_{k=0}^{\infty }$ generated by the subgradient extragradient method
(\ref{alg:cgr}) and (\ref{as}) weakly converges to a solution of the
variational inequality (\ref{a}).
\end{theorem}

\subsection{The subgradient extragradient method with outer perturbations}

\qquad In this subsection, we present the subgradient extragradient method
with outer perturbations.

\begin{algorithm}
\label{al1}$\left. {}\right. $\textbf{The subgradient extragradient method
with outer perturbations}

\textbf{Step 0}: Select a starting point $x^{0}\in\mathcal{H}$
and set $k=0$.

\textbf{Step 1}: Given the current iterate $x^{k},$ compute%
\begin{equation}
y^{k}=P_{C}(x^{k}-\gamma _{k}F(x^{k})+e_{1}(x^{k})),  \label{sStep1}
\end{equation}%
where $\gamma _{k}=\sigma \rho ^{m_{k}}$, $\sigma >0,$ $\rho \in (0,1)$ and $%
m_{k}$ is the smallest nonnegative integer such that (see \cite{Khobotov})
\begin{equation}
\gamma _{k}\Vert F(x^{k})-F(y^{k})\Vert \leq \mu \Vert x^{k}-y^{k}\Vert
,\quad \mu \in (0,1).  \label{sae}
\end{equation}%
Construct the set
\begin{equation}
T_{k}:=\{w\in \mathcal{H}|\langle (x^{k}-\gamma
_{k}F(x^{k})+e_{1}(x^{k}))-y^{k},w-y^{k}\rangle \leq 0\},  \label{g66}
\end{equation}%
and calculate
\begin{equation}
x^{k+1}=P_{T_{k}}(x^{k}-\gamma _{k}F(y^{k})+e_{2}(x^{k})).  \label{sStep2}
\end{equation}

\textbf{Step 2}:\textbf{\ }If\textbf{\ }$x^{k}=y^{k},$ then stop. Otherwise,
set $k\leftarrow (k+1)$ and return to \textbf{Step 1.}
\end{algorithm}

Denote $e_{i}^{k}:=e_{i}(x^{k})$, $i=1,2$. The sequences of perturbations $%
\{e_{i}^{k}\}_{k=0}^{\infty }$, $i=1,2$, are assume to be summable, i.e.,
\begin{equation}
\sum_{k=0}^{\infty }\Vert e_{i}^{k}\Vert <+\infty ,\quad i=1,2.  \label{BP1}
\end{equation}

Following the proof of Theorems \ref{Th31} and \ref{Th32}, we get the
convergence analysis and convergence rate of Algorithm \ref{al1}.

\begin{theorem}
\label{th62} Assume that conditions \ref{con:Condition 1.1},
 \ref{con:Condition 1.4} and \ref{con:Condition 1.5} hold. Then the sequence
$\{x^{k}\}_{k=0}^{\infty }$ generated by Algorithm \ref{al1} converges
weakly to a solution of the variational inequality (\ref{a}).
\end{theorem}

\begin{theorem}
\label{Th63} Assume that Conditions \ref{con:Condition 1.1},
\ref{con:Condition 1.4} and \ref{con:Condition 1.5} hold. Let the sequences
$\{x^k\}_{k=0}^\infty$ and $\{y^k\}_{k=0}^\infty$ be generated by Algorithm \ref{al1}.
 For any integer $t > 0$, we have a $y_t\in C$ which satisfies
\begin{equation}  \label{m-1}
\aligned \langle F(x),y_t-x\rangle\leq\frac{1}{2\Upsilon_t}%
(\|x-x^0\|^2+M(x)),\quad \forall x\in C,\endaligned
\end{equation}
where
\begin{equation}  \label{m0}
y_t=\frac{1}{\Upsilon_t}\sum_{k=0}^t\gamma_ky^k,\quad
\Upsilon_t=\sum_{k=0}^t\gamma_k,
\end{equation}
and
\begin{equation}
 M(x)=\sup_k\{\max\{\|x^{k+1}-y^k\|,\|x^{k+1}-x\|\}\}\sum_{k=0}^{\infty}\Vert e^{k}\Vert.
\end{equation}
\end{theorem}

\subsection{The BPR of the subgradient extragradient method}

\qquad In this subsection, we investigate the bounded perturbation
resilience of the subgradient extragradient method (\ref{alg:cgr}).

To this end, we treat the right-hand side of (\ref{alg:cgr}) as the algorithmic
operator $\mathbf{A}_\Psi$ of Definition \ref{Def41}, namely, we define for
all $k\geq0,$
\begin{equation}  \label{g61}
\mathbf{A}_\Psi(x^k)=P_{T(x^k)}(x^k-\gamma_kF(P_C(x^k-\gamma_k F(x^k)))),
\end{equation}
where $\gamma_k$ satisfies (\ref{as}) and
\begin{equation}  \label{g62}
T(x^k)=\{w\in \mathcal{H}|\langle(x^k-\gamma_k F(x^k))-y^k, w-y^k\rangle\leq0\}.
\end{equation}
Identify the solution set $\Psi$ with the solution set of the variational
inequality (\ref{a}) and identify the additional set $\Theta$ with $C$.

According to Definition \ref{Def41}, we need to show the convergence of the
sequence $\{x^{k}\}_{k=0}^{\infty }$ that, starting from any $x^{0}\in \mathcal{H}$,
is generated by
\begin{equation}
x^{k+1}=P_{T(x^{k}+\lambda _{k}v^{k})}((x^{k}+\lambda _{k}v^{k})-\gamma
_{k}F(P_{C}((x^{k}+\lambda _{k}v^{k})-\gamma _{k}F(x^{k}+\lambda
_{k}v^{k})))),  \label{g63}
\end{equation}%
which can be rewritten as%
\begin{equation}
\label{mEGs}
\left\{
\aligned
&y^{k} =P_{C}((x^{k}+\lambda _{k}v^{k})-\gamma _{k}F((x^{k}+\lambda
_{k}v^{k}))\medskip   \\
&T(x^{k}+\lambda _{k}v^{k}) =\{w\in \mathcal{H}|\langle ((x^{k}+\lambda
_{k}v^{k})-\gamma _{k}F(x^{k}+\lambda _{k}v^{k}))-y^{k},\\
&\qquad \qquad \qquad \quad w-y^{k}\rangle \leq
0\}\medskip    \\
&x^{k+1} =P_{T(x^{k}+\lambda _{k}v^{k})}((x^{k}+\lambda _{k}v^{k})-\gamma
_{k}F(y^{k}))
\endaligned
\right.
\end{equation}
where $\gamma _{k}=\sigma \rho ^{m_{k}}$, $\sigma >0,$ $\rho \in (0,1)$ and $%
m_{k}$ is the smallest nonnegative integer such that
\begin{equation}
\gamma _{k}\Vert F(x^{k}+\lambda _{k}v^{k})-F(y^{k})\Vert \leq \mu (\Vert
x^{k}-y^{k}\Vert +\lambda _{k}\Vert v^{k}\Vert ),\quad \mu \in (0,1).
\label{sa1s}
\end{equation}%
The sequences $\{v^{k}\}_{k=0}^{\infty }$ and $\{\lambda _{k}\}_{k=0}^{\infty }$ obey
the conditions (i) and (ii) in Definition \ref{Def41}, respectively, and
also (iii) in Definition \ref{Def41} is satisfied.\vskip0.1cm

The next theorem establishes the bounded perturbation resilience of the
subgradient extragradient method.  Since its proof is similar with that of Theorem %
\ref{Th41}, we omit it.

\begin{theorem}
\label{Th63} Assume that Conditions \ref{con:Condition 1.1}, \ref{con:Condition 1.4} and \ref{con:Condition 1.5} hold. Assume the
sequence $\{v^k\}_{k=0}^\infty$ is bounded, and the scalars $%
\{\lambda_k\}_{k=0}^\infty$ are such that $\lambda_k\geq 0$ for all $k \geq
0 $, and $\sum_{k=0}^\infty \lambda_k<+\infty$. Then the sequence $\{x^k\}_{k=0}^\infty$
generated by (\ref{mEGs}) and (\ref{sa1s}) converges weakly to a solution of
the variational inequality $(\ref{a})$.
\end{theorem}

We also get the convergence rate of the subgradient extragradient methods
with BP (\ref{mEGs}) and (\ref{sa1s}).

\begin{theorem}
\label{Th64} Assume that Conditions \ref{con:Condition 1.1}, \ref{con:Condition 1.4} and \ref{con:Condition 1.5} hold. Assume the
sequence $\{v^k\}_{k=0}^\infty$ is bounded, and the scalars $%
\{\lambda_k\}_{k=0}^\infty$ are such that $\lambda_k\geq 0$ for all $k \geq
0 $, and $\sum_{k=0}^\infty \lambda_k<+\infty$. Let the sequences $\{x^k\}_{k=0}^\infty$
and $\{y^k\}_{k=0}^\infty$ be generated by by (\ref{mEGs}) and (\ref{sa1s}). For any
integer $t > 0$, we have a $y_t\in C$ which satisfies
\begin{equation}  \label{m-1}
\aligned \langle F(x),y_t-x\rangle\leq\frac{1}{2\Upsilon_t}%
(\|x-x^0\|^2+M(x)),\quad \forall x\in C,\endaligned
\end{equation}
where
\begin{equation}  \label{m0}
y_t=\frac{1}{\Upsilon_t}\sum_{k=0}^t\gamma_ky^k,\quad
\Upsilon_t=\sum_{k=0}^t\gamma_k,
\end{equation}
and
\begin{equation}
\label{Th42-2}
M(x)=\sup_k\{\max\{\|x^{k+1}-y^k\|,3\|x^{k+1}-x\|^2\}\}\sum_{k=0}^\infty\lambda_k\|v^k\|.
\end{equation}
\end{theorem}

\subsection{Construction of the inertial subgradient extragradient methods
by BPR}

\qquad In this subsection, we construct two class of inertial subgradient
extragradient methods by using BPR, i.e., identifying the $e_i^k$, $k=1,2$
and $\lambda_k$, $v^k$ with special values.

Based on Algorithm \ref{al1}, we construct the following inertial subgradient
extragradient method:%
\begin{equation}\label{Alg:isgem}
\left\{
\aligned
&y^{k} =P_{C}(x^{k}-\gamma _{k}F(x^{k})+\alpha _{k}^{(1)}(x^{k}-x^{k-1}))) \\
&T_{k} :=\{w\in \mathcal{H}|\langle (x^{k}-\gamma _{k}F(x^{k})+\alpha
_{k}^{(1)}(x^{k}-x^{k-1}))-y^{k},\\
&\qquad \quad w-y^{k}\rangle \leq 0\},   \\
&x^{k+1}=P_{T_{k}}(x^{k}-\gamma _{k}F(y^{k})+\alpha
_{k}^{(2)}(x^{k}-x^{k-1})),
\endaligned
\right.
\end{equation}
where $\gamma _{k}$ satisfies (\ref{sa1s}) and
\begin{equation}
\alpha _{k}^{(i)}=\left\{
\begin{array}{ll}
\frac{\beta _{k}^{(i)}}{\left\Vert x^{k}-x^{k-1}\right\Vert }, & \text{if }%
\left\Vert x^{k}-x^{k-1}\right\Vert >1,\text{ }i=1,2 \\
\beta _{k}^{(i)}, & \text{if }\left\Vert x^{k}-x^{k-1}\right\Vert \leq 1.%
\end{array}%
\right.
\end{equation}

Similarly with the proof of Theorem \ref{Th51}, we get the convergence of the inertial subgradient
extragradient method (\ref{Alg:isgem}).

\begin{theorem}
\label{Th51} Assume that Conditions \ref{con:Condition 1.1}, \ref{con:Condition 1.4} and \ref{con:Condition 1.5} hold. Assume that the
sequences $\{\beta _{k}^{(i)}\}_{k=0}^{\infty }$, $i=1,2$ satisfy $%
\sum_{k=1}^{\infty }\beta _{k}^{(i)}<\infty $, $i=1,2$. Then the sequence $%
\{x^{k}\}_{k=0}^{\infty }$ generated by the inertial subgradient
extragradient method (\ref{Alg:isgem}) converges weakly to a solution of the variational
inequality $(\ref{a})$.
\end{theorem}

Using the subgradient extragradient method with bounded perturbations (\ref%
{mEGs}), we construct the following inertial subgradient extragradient
method:%
\begin{equation}\label{Alg:II-isgem}
\left\{
\aligned
&w^{k} =x^{k}+\alpha _{k}(x^{k}-x^{k-1}), \\
&y^{k} =P_{C}(w^{k}-\gamma _{k}F(w^{k})),   \\
&T_{k} :=\{w\in \mathcal{H}|\langle (w^{k}-\gamma _{k}F(w^{k}))-y^{k},w-y^{k}\rangle
\leq 0\}   \\
&x^{k+1}=P_{T_{k}}(w^{k}-\gamma _{k}F(y^{k}))
\endaligned
\right.
\end{equation}
where $\gamma _{k}=\sigma \rho ^{m_{k}}$, $\sigma >0,$ $\rho \in (0,1)$ and $%
m_{k}$ is the smallest nonnegative integer such that
\begin{equation}
\gamma _{k}\Vert F(w^{k})-F(y^{k})\Vert \leq \mu \Vert w^{k}-y^{k}\Vert
,\quad \mu \in (0,1),  \label{sa1z}
\end{equation}%
and
\begin{equation}
\alpha _{k}=\left\{
\begin{array}{ll}
\frac{\beta _{k}}{\left\Vert x^{k}-x^{k-1}\right\Vert }, & \text{if }%
\left\Vert x^{k}-x^{k-1}\right\Vert >1,\text{ }i=1,2 \\
\beta _{k}, & \text{if }\left\Vert x^{k}-x^{k-1}\right\Vert \leq 1.%
\end{array}%
\right.
\end{equation}%
We extend Theorem \ref{Th41} to the convergence of the inertial subgradient
extragradient method (\ref{Alg:II-isgem}).

\begin{theorem}
\label{Th52} Assume that Conditions \ref{con:Condition 1.1}, \ref{con:Condition 1.4} and \ref{con:Condition 1.5} hold. Assume that the
sequence $\{\beta _{k}\}_{k=0}^{\infty }$ satisfies $\sum_{k=1}^{\infty
}\beta _{k}<\infty $. Then the sequence $\{x^{k}\}_{k=0}^{\infty }$
generated by the inertial subgradient extragradient method (\ref{Alg:II-isgem}) converges
weakly to a solution of the variational inequality (\ref{a}).
\end{theorem}

\section{Numerical experiments\label{sec:num}}

\qquad In this section, we provide three examples to compare the inertial
extragradient method (\ref{Alg:iegm}) (iEG1), the inertial extragradient method (\ref{Alg:I-iegm})
(iEG2), the inertial extragradient method (\ref{Alg:II-iegm}) (iEG), the extragradient method (\ref{EG-1}), the inertial subgradient extragradient method (\ref{Alg:isgem})
(iSEG1), the inertial subgradient extragradient method (\ref{Alg:II-isgem}) (iSEG2) and the
subgradient extragradient method (\ref{alg:cgr}). \vskip 2mm

%% the following is obtained from yuchao tang, please delete. Thank you!

\begin{figure}[H]
\setlength{\floatsep}{0pt} \setlength{\abovecaptionskip}{-10pt}
\centering
\scalebox{0.6} {\includegraphics{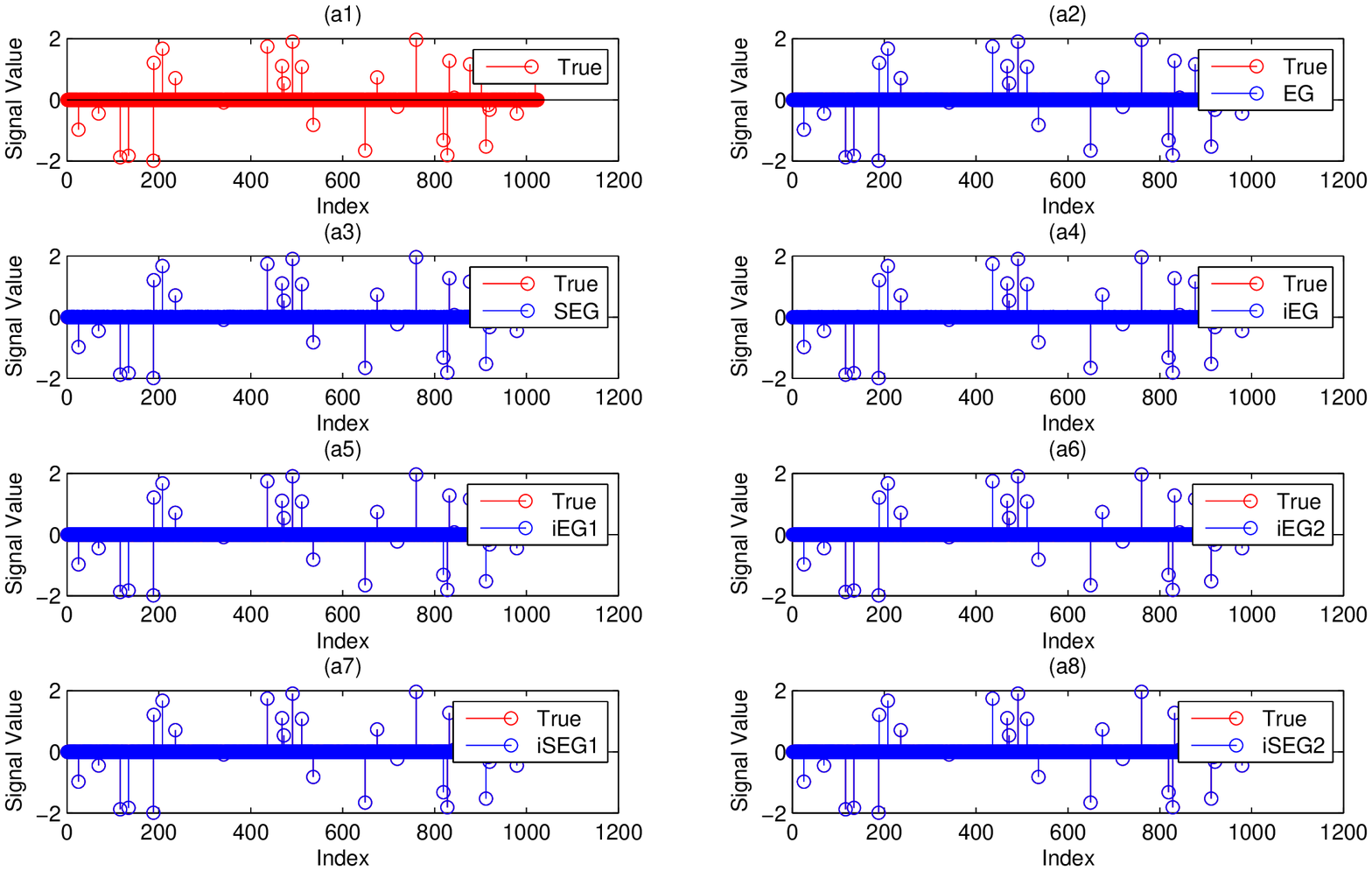}}
\caption{(a1) is the true sparse signal, (a2)-(a8) are the recovered signal vs the true signal by "EG", "SEG", "iEG", "iEG1", "iEG2" "iSEG1" and "iSEG2", respectively. }\label{ex1_k30}
\end{figure}

\begin{figure}[H]
\setlength{\floatsep}{0pt} \setlength{\abovecaptionskip}{-5pt}
\centering
\scalebox{0.6} {\includegraphics{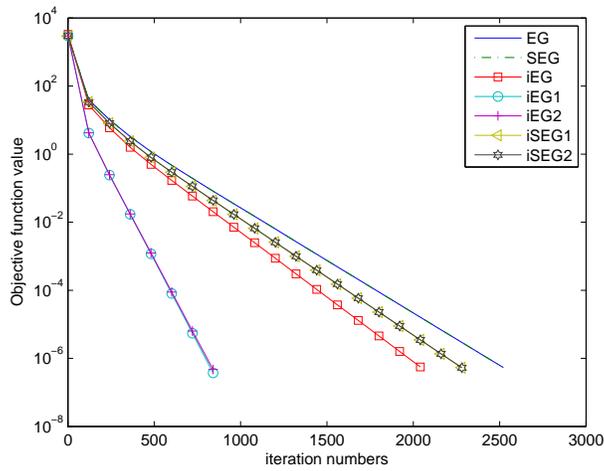}}
\caption{Comparison of the objective function value versus the iteration numbers of different methods. }\label{ex1_k30_fun}
\end{figure}

In the first example, we consider a typical sparse signal recovery problem. We choose the following set of parameters. Take $\sigma =5$, $\rho =0.9$ and $\mu =0.7.$ Set
\begin{equation}
\alpha _{k}=\alpha _{k}^{(i)}=  \label{fs2} \\
\frac{1}{k^{2}}\quad \hbox{if}\quad \Vert x^{k}-x^{k-1}\Vert \leq 1,
\end{equation}%
in inertial extragradient methods (\ref{Alg:iegm}) and (\ref{Alg:I-iegm}), and inertial subgradient
extragradient methods (\ref{Alg:isgem}) and (\ref{Alg:II-isgem}). Choose $\alpha _{k}=0.35$ and $\lambda
_{k}=0.8$ in the inertial extragradient method (\ref{Alg:I-iegm}).
\begin{example}\label{ex:sparse}
Let $x_0 \in R^{n}$ be a $K$-sparse signal, $K<<n$. The sampling matrix $A \in R^{m\times n} (m<n)$ is stimulated by standard Gaussian distribution
 and vector $b = Ax_0 + e$, where $e$ is additive noise. When $e=0$, it means that there is no noise to the observed data. Our task is to recover the signal $x_0$ from the data $b$.
\end{example}
It's well-known that the sparse signal $x_0$ can be recovered by solving the following LASSO problem \cite{Tibshirani96},
\begin{equation}
\begin{aligned}
\min_{x\in R^n}\ & \frac{1}{2}\| Ax - b \|_{2}^{2} \\
s.t. \ & \|x\|_1 \leq t,
\end{aligned}
\end{equation}
where $t >0$. It is easy to see that the optimization problem (\ref{ex:sparse}) is a special case of the variational inequality problem (\ref{a}), where $F(x) = A^{T}(Ax-b)$ and $C = \{ x | \|x\|_1 \leq t \}$. We can use the proposed iterative algorithms to solve the optimization problem (\ref{ex:sparse}). Although the orthogonal projection onto the closed convex set $C$ doesn't have a closed-form solution, the projection operator $P_{C}$ can be precisely computed in a polynomial time. We include the detail of computing $P_{C}$ in the Appendix.
We conduct plenty of simulations to compare the performance of the proposed iterative algorithms. The following inequality was defined as the stopping criteria,
$$
\| x^{k+1} - x^k \| \leq \epsilon,
$$
where $\epsilon >0$ is a given small constant. $"Iter"$ denotes the iteration numbers. $"Obj"$ represents the objective function value and $"Err"$ is the $2$-norm error between the recovered signal and the true $K$-sparse signal. We divide the experiments into two parts. One task is to recover the sparse signal $x_0$ from noise observation vector $b$ and the other is to recover the sparse signal from noiseless data $b$. For the noiseless case, the obtained numerical results are reported in Table \ref{example1:results1}. To visually view the results, Figure \ref{ex1_k30} shows the recovered signal compared with the true signal $x_0$ when $K=30$. We can see from Figure \ref{ex1_k30} that the recovered signal is the same as the true signal. Further, Figure \ref{ex1_k30_fun} presents the objective function value versus the iteration numbers.

\begin{table}[H]
\footnotesize
\centering
\caption{Numerical results obtained by the proposed iterative algorithms
when $m =240, n =1024$ in the noiseless case.}
\begin{tabular}{c|c|cccccccc}
\hline
\multirow{2}[1]{*}{$K$-sparse}  & \multirow{2}[1]{*}{Methods} & &   \multicolumn{3}{c}{$ \epsilon = 10^{-4}$} &  & \multicolumn{3}{c}{$ \epsilon = 10^{-6}$} \\
\cline{4-6} \cline{8-10}
signal & & & $Iter$ & $Obj$ & $Err$ & &  $Iter$ &  $Obj$ & $Err$  \\
\hline
 \hline
\multirow{7}[1]{*}{$K=20$}  & EG & & $444$ & $9.7346e-4$ & $0.0080$ & & $817$ & $9.6625e-8$ & $7.9856e-5$\\
& SEG & & $444$ & $9.7272e-4$ & $0.0080$ & & $817$ & $9.6555e-8$ & $7.9827e-5$ \\
& iEG & & $374$ & $6.2389e-4$ & $0.0064$ & & $675$ & $6.3456e-8$ & $6.4715e-5$ \\
& iEG1 & & $159$ & $7.0799e-5$ & $0.0021$ & & $263$ & $7.4280e-9$ & $2.2041e-5$ \\
& iEG2 & & $158$ & $8.3897e-5$ & $0.0023$ & & $273$ & $1.0889e-8$ & $2.6809e-5$ \\
& iSEG1 & & $415$ & $8.9563e-4$ & $0.0076$ & & $787$ & $2.3571e-7$ & $5.2470e-5$ \\
& iSEG2 & & $414$ & $9.2167e-4$ & $0.0077$ & & $760$ & $9.1586e-8$ & $7.7275e-5$ \\
\hline
\multirow{7}[1]{*}{$K=30$}  & EG & & $1285$ & $0.0035$ & $0.0281$ & & $2583$ & $3.4535e-7$ & $2.8035e-4$\\
& SEG & & $1285$ & $0.0035$ & $0.0281$ & & $2583$ & $3.4534e-7$ & $2.8035e-4$ \\
& iEG & & $1091$ & $0.0023$ & $0.0227$ & & $2144$ & $2.2732e-7$ & $2.2745e-4$ \\
& iEG1 & & $532$ & $3.7493e-4$ & $0.0092$ & & $944$ & $3.7522e-8$ & $9.2287e-5$ \\
& iEG2 & & $535$ & $3.7961e-4$ & $0.0093$ & & $956$ & $4.3181e-8$ & $9.3120e-5$ \\
& iSEG1 & & $1176$ & $0.0031$ & $0.0266$ & & $2351$ & $3.1038e-7$ & $2.6137e-4$ \\
& iSEG2 & & $1176$ & $0.0031$ & $0.0266$ & & $2346$ & $3.1635e-7$ & $2.6784e-4$ \\
\hline
\multirow{7}[1]{*}{$K=40$}  & EG & & $1729$ & $0.0050$ & $0.0405$ & & $3599$ & $5.0237e-7$ & $4.0488e-4$\\
& SEG & & $1729$ & $0.0050$ & $0.0405$ & & $3599$ & $5.0228e-7$ & $4.0484e-4$ \\
& iEG & & $1473$ & $0.0033$ & $0.0328$ & & $2990$ & $3.3182e-7$ & $3.2905e-4$ \\
& iEG1 & & $744$ & $5.4838e-4$ & $0.0134$ & & $1361$ & $5.5456e-8$ & $1.3440e-4$ \\
& iEG2 & & $745$ & $5.4807e-4$ & $0.0134$ & & $1355$ & $6.4785e-8$ & $1.4191e-4$ \\
& iSEG1 & & $1570$ & $0.0045$ & $0.0384$ & & $3246$ & $4.5079e-7$ & $3.8146e-4$ \\
& iSEG2 & & $1572$ & $0.0045$ & $0.0382$ & & $3244$ & $4.5389e-7$ & $3.8435e-4$ \\

 \hline
\hline

\end{tabular}\label{example1:results1}
\end{table}

\begin{figure}[H]
\setlength{\floatsep}{0pt} \setlength{\abovecaptionskip}{-5pt}
\centering
\scalebox{0.6} {\includegraphics{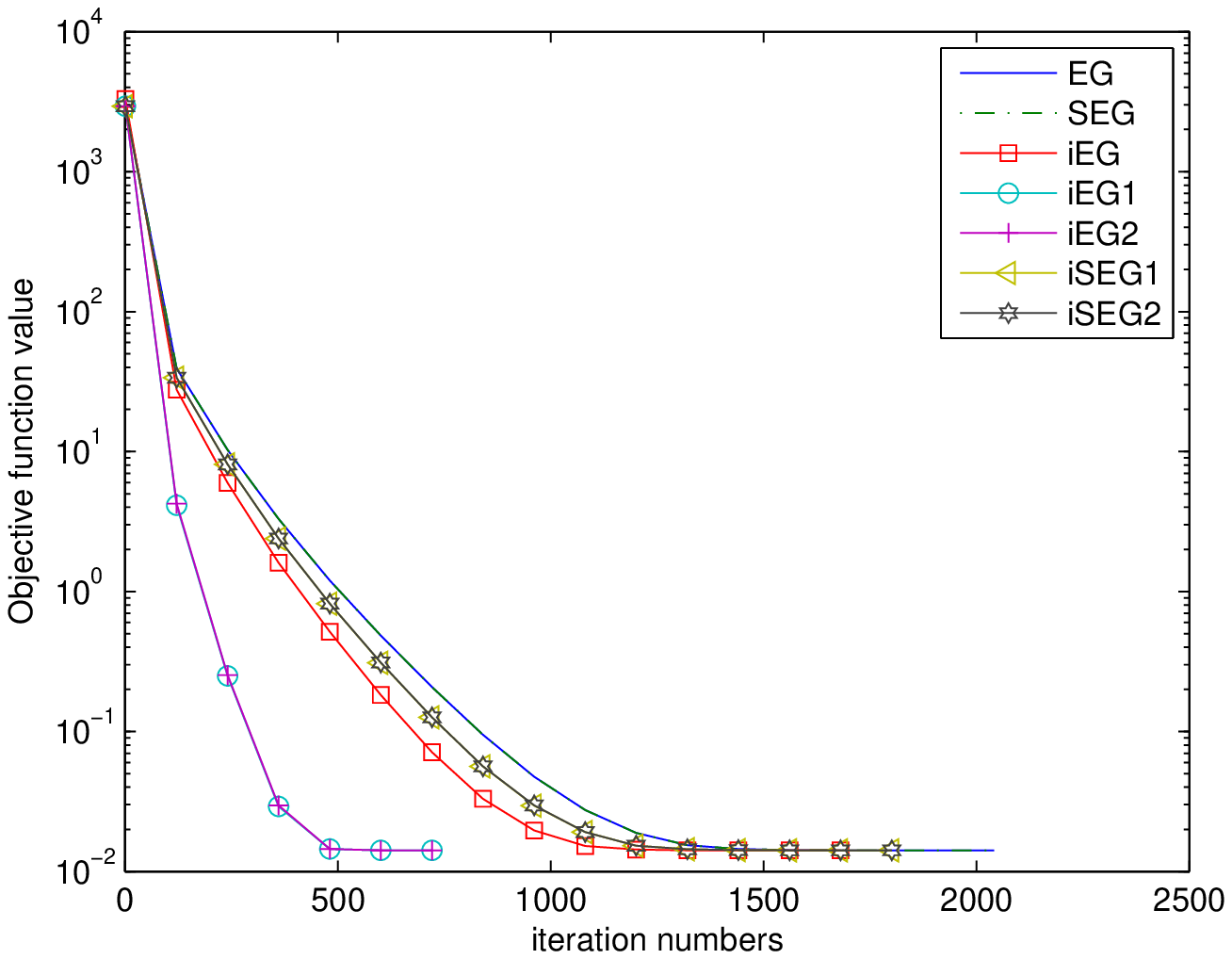}}
\caption{Comparison of the objective function value versus the iteration numbers of different methods in the noise case of $\beta = 0.02$. }\label{ex1_k30_fun_noise}
\end{figure}

\begin{figure}[H]
\setlength{\floatsep}{0pt} \setlength{\abovecaptionskip}{-10pt}
\centering
\scalebox{0.6} {\includegraphics{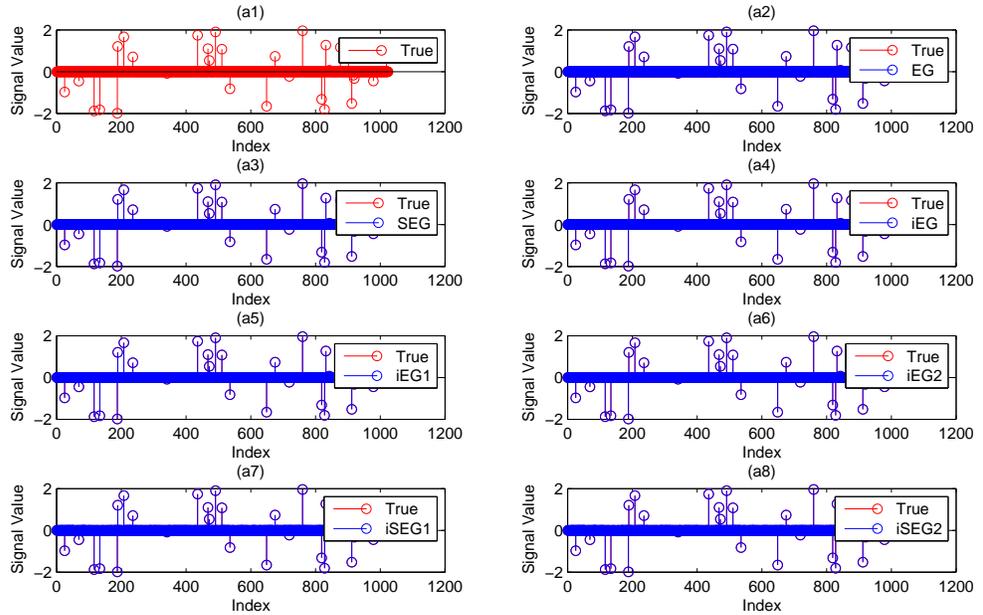}}
\caption{(a1) is the true sparse signal, (a2)-(a8) are the recovered signal vs the true signal by "EG", "SEG", "iEG", "iEG1", "iEG2" "iSEG1" and "iSEG2" in the noise case of $\beta = 0.02$, respectively. }\label{ex1_k30_noise}
\end{figure}

For the noise observation $b$, we assume that the vector $e$ is corrupted by Gaussian noise with zero mean and $\beta$ variances. The system matrix $A$ is the same as the noiseless case and the sparsity level $K=30$.
We list the numerical results for different noise level $\beta$ in Table \ref{example1:results2}. When the noise $\beta = 0.02$,  Figure \ref{ex1_k30_fun_noise} shows the objective function value versus the iteration numbers.  Figure \ref{ex1_k30_noise} shows the recovered signal vs the true signal in the noise case.

\begin{table}[H]
\footnotesize
\centering
\caption{Numerical results for the proposed iterative algorithms with different noise value $\beta$.}
\begin{tabular}{c|c|cccccccc}
\hline
\multirow{2}[1]{*}{Variances}  & \multirow{2}[1]{*}{Methods} & &   \multicolumn{3}{c}{$ \epsilon = 10^{-4}$} &  & \multicolumn{3}{c}{$ \epsilon = 10^{-6}$} \\ \cline{4-6} \cline{8-10}
$\beta = $ & & & $Iter$ & $Obj$ & $Err$ & &  $Iter$ &  $Obj$ & $Err$  \\
\hline
 \hline
\multirow{7}[1]{*}{$0.01$}  & EG & & $1264$ & $0.0092$ & $0.0317$ & & $2192$ & $0.0061$ & $0.0131$\\
& SEG & & $1264$ & $0.0092$ & $0.0317$ & & $2192$ & $0.0061$ & $0.0131$ \\
& iEG & & $1070$ & $0.0081$ & $0.0272$ & & $1812$ & $0.0061$ & $0.0131$ \\
& iEG1 & & $519$ & $0.0063$ & $0.0164$ & & $788$ & $0.0061$ & $0.0130$ \\
& iEG2 & & $516$ & $0.0063$ & $0.0166$ & & $786$ & $0.0061$ & $0.0130$ \\
& iSEG1 & & $1156$ & $0.0089$ & $0.0305$ & & $1995$ & $0.0061$ & $0.0131$ \\
& iSEG2 & & $1157$ & $0.0089$ & $0.0304$ & & $1990$ & $0.0061$ & $0.0131$ \\
\hline
\multirow{7}[1]{*}{$0.02$}  & EG & & $1274$ & $0.0163$ & $0.0387$ & & $2086$ & $0.0142$ & $0.0272$\\
& SEG & & $1274$ & $0.0163$ & $0.0387$ & & $2086$ & $0.0142$ & $0.0272$\\
& iEG & & $1070$ & $0.0154$ & $0.0356$ & & $1728$ & $0.0142$ & $0.0272$\\
& iEG1 & & $492$ & $0.0144$ & $0.0300$ & & $756$ & $0.0142$ & $0.0272$\\
& iEG2 & & $495$ & $0.0143$ & $0.0300$ & & $759$ & $0.0142$ & $0.0272$\\
& iSEG1 & & $1163$ & $0.0161$ & $0.0378$ & & $1899$ & $0.0142$ & $0.0272$\\
& iSEG2 & & $1161$ & $0.0161$ & $0.0380$ & & $1895$ & $0.0142$ & $0.0272$\\
\hline
\multirow{7}[1]{*}{$0.05$}  & EG & & $1190$ & $0.1012$ & $0.0749$ & & $1869$ & $0.0991$ & $0.0651$\\
& SEG & & $1190$ & $0.1012$ & $0.0749$ & & $1869$ & $0.0991$ & $0.0651$ \\
& iEG & & $996$ & $0.1005$ & $0.0727$ & & $1542$ & $0.0991$ & $0.0650$ \\
& iEG1 & & $460$ & $0.0993$ & $0.0677$ & & $670$ & $0.0991$ & $0.0650$ \\
& iEG2 & & $461$ & $0.0993$ & $0.0675$ & & $665$ & $0.0991$ & $0.0650$ \\
& iSEG1 & & $1084$ & $0.1010$ & $0.0742$ & & $1704$ & $0.0991$ & $0.0651$ \\
& iSEG2 & & $1084$ & $0.1010$ & $0.0742$ & & $1704$ & $0.0991$ & $0.0651$ \\

 \hline
\hline

\end{tabular}\label{example1:results2}
\end{table}

\begin{example} \label{ex:Example 6.1}, Let $F:\mathbb{R}^2 \rightarrow \mathbb{R}^2$ be
defined by
\begin{equation}  \label{fs1}
F(x,y)=(2x+2y+\sin(x),-2x+2y+\sin(y)),\quad \forall x,y \in \mathbb{R}.
\end{equation}
\end{example}

The authors \cite{Dong1} proved that $F$ is Lipschitz continuous with $L=%
\sqrt{26}$ and 1-strongly monotone. Therefore the variational inequality (%
\ref{a}) has a unique solution and $(0,0)$ is its solution.

Let $C=\{x\in \mathbb{R}^2\,|\,e_1\leq x\leq e_2\}$, where $e_1=(-10,-10)$
and $e_2=(100,100)$. Take the initial point $x_0=(-100,10)\in \mathbb{R}^2$.
Since $(0,0)$ is the unique solution of the variational inequality (\ref{a}%
), denote by $D_k:=\|x^k\|\leq 10^{-5}$ the stopping criterion.

\begin{figure}[H]
\begin{center}
\includegraphics[width=10.2cm,height=7.2cm]{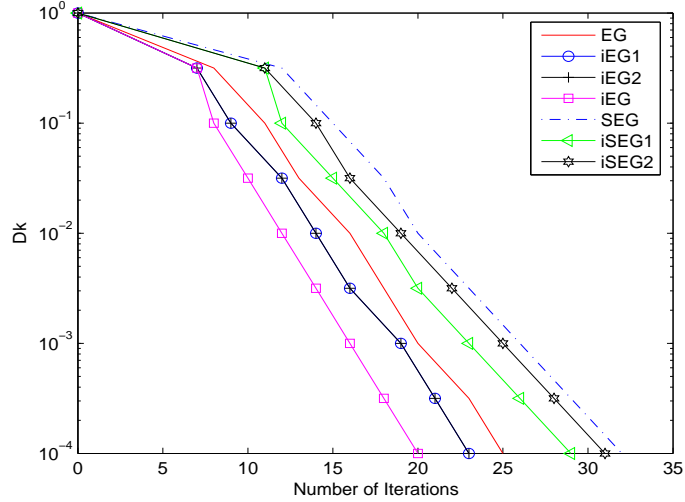}
\end{center}
\caption{Comparison of the number of iterations of different methods for
Example \ref{ex:Example 6.1}}\label{Fig:1}
\end{figure}

\begin{example} \label{ex:Example 6.2}
Let $F: \mathbb{R}^n \rightarrow \mathbb{R}^n$
defined by $F(x)=Ax+b$, where $A=Z^TZ$, $Z=(z_{ij})_{n\times n}$ and $%
b=(b_i)\in\mathbb{\ R}^n$ where $z_{ij}\in (0,1)$ and $b_i\in(0,1)$ are
generated randomly.
\end{example}

It is easy to verify that $F$ is $L-$Lipschitz
continuous and $\eta-$strongly monotone with $L=\max(eig(A))$ and $%
\eta=\min(eig(A))$.

Let $C:=\{x\in \mathbb{R}^{n}\,|\,\Vert x-d\Vert \leq r\}$, where the center
\begin{equation}
d\in \lbrack (-10,-10,\ldots ,-10),(10,10,\ldots ,10)]\subset \mathbb{R}^{n}
\end{equation}
and radius $r\in (0,10)$ are randomly chosen. Take the initial point $%
x_{0}=(c_{i})\in \mathbb{R}^{n}$, where $c_{i}\in \lbrack 0,2]$ is generated
randomly. Set $n=100$. Take $\rho =0.4$ and other parameters are set the
same values as Example \ref{ex:Example 6.1}. Although the variational inequality (\ref{a})
has an unique solution, it is difficult to get the exact solution. So,
denote by $D_{k}:=\Vert x^{k+1}-x^{k}\Vert \leq 10^{-5}$ the stopping
criterion.

\begin{figure}[H]
\begin{center}
\includegraphics[width=10.2cm,height=7.2cm]{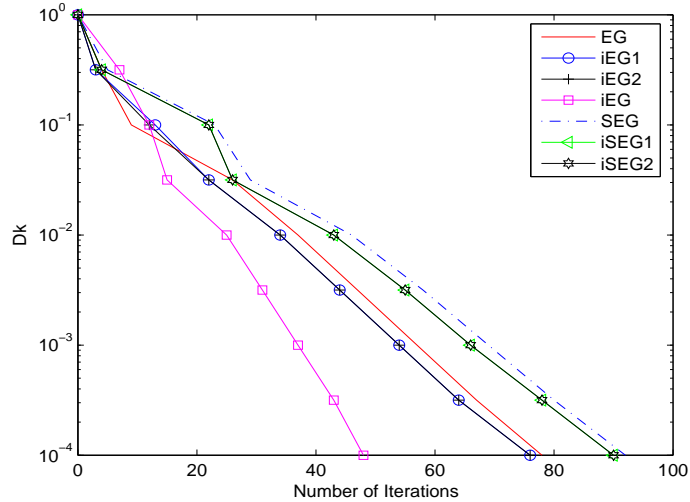}
\end{center}
\caption{Comparison of the number of iterations of different methods for
Example \ref{ex:Example 6.2}}\label{Fig:2}
\end{figure}

From Figures \ref{Fig:1} and \ref{Fig:2}, we conclude: (i) The inertial type algorithms improves
the original algorithms; (ii) the performance of the inertial extragradient
methods (\ref{Alg:iegm}) and (\ref{Alg:I-iegm}) are almost the same; (iii) the inertial subgradient
extragradient method (\ref{Alg:isgem}) performs better than the inertial subgradient
extragradient method (\ref{Alg:II-isgem}) for Example 6.1, while they are almost the same
for Example 6.2; (iv) the (inertial) extragradient methods behave better
than the (inertial) subgradient extragradient methods since the sets $C$ in
Examples \ref{ex:Example 6.1} and \ref{ex:Example 6.2} are simple and hence the computation load of the
projection onto it is small; (v) the inertial extragradient method (\ref{Alg:iegm}) has an advantage over the inertial extragradient methods (\ref{Alg:iegm}) and (\ref{Alg:I-iegm}). The reason may be that it takes bigger the inertial parameter $\alpha _{k}.$

\section*{Appendix}

In this part, we present the detail of computing a vector $y\in R^n$ onto the $\ell_1$-norm ball constraint. For convenience, we consider projection onto the unit $\ell_1$-norm ball first. Then we extend it to the general $\ell_1$-norm ball constraint.

The projection onto the unit $\ell_1$-norm ball is to solve the optimization problem,
\begin{equation*}\label{l1-ball}
\begin{aligned}
\min_{x\in R^n}\ & \frac{1}{2}\| x- y \|_{2}^{2} \\
s.t. \ & \|x\|_1 \leq 1.
\end{aligned}
\end{equation*}
The above optimization problem is a typical constrained optimization problem, we consider to solve it based on the Lagrangian method.
Define the Lagrangian function $L(x,\lambda)$ as
\begin{equation*}
L(x,\lambda) = \frac{1}{2}\|  x -y  \|_{2}^{2} + \lambda (  \|x\|_1 - 1 ).
\end{equation*}
Let $(x^*, \lambda^*)$ be the optimal primal and dual pair. It satisfies the KKT conditions of
\begin{equation*}
\begin{aligned}
& 0 \in  (x^* - y) + \lambda^* \partial(\|x^*\|_1) \\
& \lambda^* ( \|x^*\|_1 - 1 ) = 0 \\
& \lambda^* \geq 0.
\end{aligned}
\end{equation*}
It is easy to check that if $\|y\|_1 \leq 1$, then $x^* = y$ and $\lambda^* =0$. In the following, we assume $\|y\|_1 >1$. Based on the KKT conditions, we obtain $\lambda^* >0$ and $\|x^*\|_1 = 1$. From the first order optimality, we have $x^{*} = \max\{ |y|-\lambda^*, 0 \}\otimes Sign(y)$, where $\otimes$ represents element-wise multiplication and $Sign(\cdot)$ denotes the symbol function, i.e., $Sign(y_i) = 1$ if $y_i \geq 0$; otherwise $Sign(y_i) = -1$.

Define a function $f(\lambda) = \| x(\lambda) \|_1$, where $x(\lambda) = S_{\lambda}(y) = \max \{ |y|-\lambda, 0 \}\otimes Sign(y)$. We prove the following lemma.

\begin{lemma}
For the function $f(\lambda)$, there must exist a $\lambda^* >0$ such that $f(\lambda^*) = 1$.
\end{lemma}

\textit{Proof.}
Since $f(0) = \| S_{0}(y) \|_1 = \|y\|_1 >1$. Let $\lambda^+ = \max_{1\leq i \leq n}\{ |y_i| \}$, then $f(\lambda^+) =0 <1$. Notice that $f(\lambda)$ is decreasing and convex. Therefore, by the intermediate value theorem, there exists $\lambda^* >0$ such that $f(\lambda^*) = 1$. $\Box$

To find a $\lambda^*$ such that $f(\lambda^*) =1$. We follow the following steps:

\noindent \textbf{Step 1.} Define a vector $\overline{y}$ with the same element as $|y|$, which was sorted in descending order. That is $\overline{y}_{1} \geq \overline{y}_2 \geq \cdots \overline{y}_{n} \geq 0$.

\noindent \textbf{Step 2.} For every $k= 1,2, \cdots, n$, solve the equation $\sum_{i=1}^{k}\overline{y}_i - k \lambda =1$. Stop search until the solution $\lambda^*$ belongs to the interval $[\overline{y}_{k+1}, \overline{y}_{k}]$.

In conclusion, the optimal $x^*$ can be computed by $x^{*} = \max\{ |y|-\lambda^*, 0 \}\otimes Sign(y)$. The next lemma extend the projection onto the unit $\ell_1$-norm ball to the general $\ell_1$-norm ball constraint.

\begin{lemma}
Let $C_1 = \{ x | \|x\|_1 \leq 1 \}$. For any $t>0$, define a general $\ell_1$-norm ball constraint set $C = \{ x | \|x \|_1 \leq t \}$. Then for any vector $y\in R^n$, we have
$$
P_{C}(y) = t P_{C_1}(\frac{y}{t}).
$$
\end{lemma}

\textit{Proof.}
To compute the projection $P_{C}(y)$, it is to solve the optimization problem,
\begin{equation*}
\begin{aligned}
P_{C}(y) = \arg \min_{x\in R^n}\ & \frac{1}{2}\|  x-y \|_{2}^{2} \\
s.t. \ & \| x \|_1 \leq t.
\end{aligned}
\end{equation*}
For any $x\in C$, let $\overline{x} = \frac{x}{t}$, it follows that $\overline{x}\in C_1$. The optimal solution $x^{*}$ of the above optimization problem satisfying $x^{*} = P_{C}(y) = t \overline{x}^{*}$, where $\overline{x}^{*}$ is the optimal solution of the optimization problem of,
\begin{equation*}
\begin{aligned}
\overline{x}^{*} = \arg \min_{\overline{x}\in R^n}\ & \frac{1}{2}\|  \overline{x}-\frac{y}{t} \|_{2}^{2} \\
s.t. \ & \| \overline{x} \|_1 \leq 1.
\end{aligned}
\end{equation*}
It is observed that $\overline{x}^{*}$ is exact projection onto the closed convex set $C_1$. That is $\overline{x}^{*} = P_{C_1}(\frac{y}{t})$. This completes the proof.  $\Box$

\section*{Conclusions}
In this research article we study an important property of iterative algorithms for solving variational inequality (VI) problems and it is called called bounded perturbation resilience. In particular we focus in extragradient-type methods. This enable use to develop inexact versions of the methods as well as applying the superiorizion methodology in order to obtain a "superior" solution to the original problem. In addition, some inertial extragradient methods are also derived. All the presented methods converge at the rate of $O(1/t)$ and three numerical examples
illustrate, demonstrate and compare the performances of all the algorithms.

\section*{Competing interests}
The authors declare that they have no competing interests.

\section*{Funding}
The first author is supported by National Natural Science Foundation of China (No. 61379102) and
Open Fund of Tianjin Key Lab for Advanced Signal Processing (No. 2016ASP-TJ01). The third author is supported by the EU FP7 IRSES program STREVCOMS, grant no. PIRSES-GA-2013-612669. The fourth author is supported by supported by Visiting Scholarship of Academy of Mathematics and Systems Science, Chinese Academy of Sciences (AM201622C04) and the National Natural Science Foundations of China (11401293,11661056), the Natural Science Foundations of Jiangxi Province (20151BAB211010).

\section*{Authors’ contributions}
All authors contributed equally to the writing of this paper. All authors read and approved
the final manuscript.

\section*{Acknowledgment}
We wish to thank the anonymous referees for their thorough analysis and
review, all their comments and suggestions helped tremendously in improving
the quality of this paper and made it suitable for publication.

%\nocite{*}

\end{document}